\documentclass[journal]{IAENGtran}
\ifCLASSINFOpdf
   \usepackage[pdftex]{graphicx}
   \DeclareGraphicsExtensions{.pdf,.jpeg,.png}
\else
   \usepackage[dvips]{graphicx}
   \DeclareGraphicsExtensions{.eps}
\fi

\usepackage{amssymb}
\usepackage{subfigure, graphicx}
\usepackage{latexsym}
\usepackage{graphicx}
\usepackage{amssymb, amsmath}
\usepackage{bm}

\def\G{\Gamma}
\def\calP{{\cal P}}
\def\calS{{\cal S}}
\def\scF{{\cal F}}
\def\R{\mathbb{R}}
\def\Z{\mathbb{Z}}
\def\N{\mathbb{N}}

\def\veca{\bm{a}}
\def\vecb{\bm{b}}
\def\vecf{\bm{f}}

\def\vecx{\bm{x}}

\def\vecy{\bm{y}}
\def\vecT{\bm{T}}
\def\vecN{\bm{N}}

\def\Image{\mathop\mathrm{Image}}

\def\<{\langle}
\def\>{\rangle}

\def\eps{\varepsilon}

\def\ds{\mathrm{d}s}
\def\du{\mathrm{d}u}
\def\dt{\mathrm{d}t}
\def\dnu{\mathrm{d}\nu}

\newtheorem{theorem}{Theorem}
\newtheorem{remark}{Remark}

\begin{document}
%
\title{On a Gradient Flow of Plane Curves Minimizing the Anisoperimetric Ratio}

\author{Daniel \v{S}ev\v{c}ovi\v{c} and Shigetoshi Yazaki
\thanks{Manuscript received April 27, 2013; revised June 6, 2013.  This work was supported in part by the VEGA grant 1/0747/12 (DS) and Grant-in-Aid for Scientific Research (C) 23540150 (SY)}
\thanks{Prof. D. \v{S}ev\v{c}ovi\v{c}, PhD., is with the 
Department of Applied Mathematics and Statistics, Comenius University, 842 48 Bratislava, Slovak Republic. 
e-mail: sevcovic@fmph.uniba.sk}
\thanks{Prof. S. Yazaki, PhD., is with the Department of Mathematics, School of Science and Technology, Meiji University, Kanagawa 214-8571, Japan. email: syazaki@meiji.ac.jp}}

\maketitle

\pagestyle{empty}
\thispagestyle{empty}

\begin{abstract}

We analyze a gradient flow of closed planar curves minimizing the anisoperimetric ratio. For such a flow the normal velocity is a function of the anisotropic curvature and it also depends on the total interfacial energy and enclosed area of the curve. In contrast to the gradient flow for the isoperimetric ratio, we show there exist initial curves for which the enclosed area is decreasing with respect to time. We also derive a mixed anisoperimetric inequality for the product of total interfacial energies corresponding to different anisotropy functions. Finally, we present several computational examples illustrating theoretical results. 
\end{abstract}

\begin{IAENGkeywords}
Anisoperimetric ratio, gradient geometric flows, mixed anisoperimetric ratio inequality, tangential stabilization
\end{IAENGkeywords}

\section{ Introduction}
\label{sec-1}

\IAENGPARstart{T}{he} goal of this paper is to investigate a geometric flow of closed plane curves $\Gamma^t, t\ge0$, minimizing the anisoperimetric ratio. We will show that the normal velocity $\beta$ for such a geometric flow is a function of the anisotropic curvature, the total interfacial energy and enclosed area of an evolved curve, 
\begin{equation}\label{geomrovnonloc}
\beta = \delta(\nu) k  + \scF_\Gamma,
\end{equation}
where $k$ is the curvature and $\delta(\nu)>0$ is a strictly positive coefficient depending on the tangent angle $\nu$ at a point $\vecx\in \Gamma^t$. Here $\scF_\Gamma$ is a nonlocal part of the normal velocity depending on the entire shape of the curve $\Gamma=\Gamma^t$ and the term $\delta(\nu) k$ represents the anisotropic curvature. In typical situations, the nonlocal part is a function of the enclosed area $A$ and the interfacial energy $L_\sigma = \int_\Gamma \sigma \ds$, i.e. $\scF_\Gamma = \scF(A,L_\sigma)$. 
As an example one can consider 
\[
\beta = k  - \frac{2\pi}{L},
\]
where  $L \equiv L_1$ is the length of an evolved closed curve $\Gamma$. It is well known that such a flow represents the area preserving geometric evolution of closed embedded plane curves investigated by Gage \cite{Gage1986}. Among other geometric flows with nonlocal normal velocity we mention the curvature driven length preserving flow in which $\beta = k - \frac{1}{2\pi} \int_\Gamma k^2 \ds$
studied by Ma and Zhu \cite{Ma2012} and the inverse curvature driven flow preserving the length $\beta = -k^{-1} + \frac{L}{2\pi}$ studied by Pan and Yang \cite{Pan2008}.
The isoperimetric ratio gradient flow with $\beta = k - L/(2A)$ has been proposed and investigated by Jiang and Zhu \cite{Jiang2008} for convex curves and by the authors in \cite{SY2011} for general closed Jordan curves evolving in the plane.

Recently, a classical nonlocal curvature flow preserving the enclosed area was reinvestigated by Xiao \emph{et al.} in \cite{Xiao2013}. They proved uniform upper bound and lower bound on the curvature. Furthermore, Mao \emph{et al.} \cite{Mao2012} showed that such a nonlocal flow will decrease the perimeter of the evolving curve and make the curve more and more circular during the evolution process. Applying inequalities of Andrews and Green-Osher type, Lin and Tsai \cite{Lin2012} showed that the evolving curves will converge to a round circle, provided that the curvature is a-priori bounded. However, most of those fine results for area preserving flow still have to be extended to the case of a class of non-local flows minimizing the isoperimetric and/or anisoperimetric ratio.

The main goal of this paper is twofold. First we derive the normal velocity $\beta$ corresponding to the anisoperimetric ratio gradient flow. It turns out that $\beta = k_\sigma - L_\sigma/(2A)$ where $k_\sigma$ is the anisotropic curvature, i.e. $\beta$ has the form of (\ref{geomrovnonloc}). We derive and analyze several important properties of such a geometric flow. In contrast to the isoperimetric ratio gradient flow 
(c.f. Jiang and Zhu \cite{Jiang2008}, \cite{SY2011}),
we show that the anisoperimetric ratio gradient flow may initially increase the total length and, conversely, decrease the enclosed area of evolved curves. In order to verify such striking phenomena, an accurate numerical discretization scheme for fine approximation of the geometric flow has to be proposed. This is the second principal goal of the paper. We derive a numerical scheme based on the method of flowing finite volumes with combination of asymptotically uniform tangential redistribution of grid points.
The idea of a uniform tangential redistribution has been proposed by How \emph{et al} in \cite{HouLS1994} and further analyzed by Mikula and \v{S}ev\v{c}ovi\v{c} in \cite{MikulaS2001}. The asymptotically uniform tangential redistribution has been analyzed in \cite{MikulaS2004b,MikulaS2004a}. 
The scheme is tested on the area-decrease and length-increase phenomena as well as on various other examples of evolution of initial curves having large variations in the curvature.

The paper is organized as follows. In the next section we recall the system of governing PDEs describing the evolution of all relevant geometric quantities. In section 3 we recall basic properties  the anisotropic curvature and Wulff shape. We prove an important duality identity between total interfacial energies corresponding to different anisotropies. In section 4 we investigate a gradient flow for the anisoperimetric ratio. It turns out that the flow of plane minimizing the anisoperimetric ratio has the normal velocity locally depending on the anisotropic curvature and nonlocally depending on the total interfacial energy and the enclosed area of the evolved curve. Section 5 is devoted to the proof of a mixed anisoperimetric inequality for the product of two total interfacial energies corresponding to two anisotropy functions. In section 6 we investigate properties  of the enclosed area for the anisoperimetric gradient flow. In contrast to a gradient flow for the isoperimetric ratio, we will show that there are initial convex curves for which the enclosed area is strictly decreasing. Finally, in section 7 we construct a counterexample to a comparison principle showing that there initial noninteresting curves such that they intersect each other immediately when evolved in the normal direction by the anisoperimetric ratio gradient flow. In section 8 we derive a numerical scheme for solving curvature driven flows with normal velocity depending on no-local terms. The scheme is based on a flowing finite volume method combined with a precise scheme for approximation of non-local terms. We present several numerical examples illustrating theoretical results and interesting phenomena for the gradient flow for anisoperimetric ratio.

\section{System of governing equations and curvature adjusted tangential redistribution}
\label{sec:GE}

In this section we recall description and basic properties of geometric evolution of a closed plane Jordan curve $\G$ which can be parameterized by a smooth function $\vecx:\ [0, 1]\to\R^2$ such that $\G=\Image(\vecx)=\{\vecx(u);\ u\in[0, 1]\}$ and  $|\partial_u\vecx|>0$. We identify the interval $[0,1]$ with the quotient space $\R/\Z$ by imposing periodic boundary conditions for $\vecx(u)$ at $u=0,1$. We denote $\partial_\xi{\sf F}=\partial{\sf F}/\partial\xi$, and 
$|\veca|=\sqrt{\veca \cdot \veca}$ where  $\veca \cdot \vecb$ is the Euclidean inner product between vectors $\veca$ and $\vecb$. The unit tangent vector is given by  $\vecT=\partial_u\vecx/|\partial_u\vecx|=\partial_s\vecx$, where $s$ is the arc-length parameter $\ds=|\partial_u\vecx|\du$. The unit inward normal vector is defined in such a way that $\det(\vecT, \vecN)=1$. Then the  signed curvature $k$ in the direction $\vecN$ is given by $k=\det(\partial_s\vecx, \partial_s^2\vecx)$. Let $\nu$ be a tangent angle, i.e., $\vecT=(\cos\nu, \sin\nu)^{\mathrm{T}}$ and $\vecN=(-\sin\nu, \cos\nu)^{\mathrm{T}}$. From the Fren\'et formulae $\partial_s \vecT = k \vecN$ and $\partial_s \vecN = -k \vecT$ we deduce that  $\partial_s\nu =k$. 

Geometric evolution problem can be formulated as follows: for a given initial curve $\G^0=\Image(\vecx^0)=\G$, find a family of curve $\{\G^t\}_{t\geq 0}$, $\G^t=\{\vecx(u, t);\ u\in[0, 1]\}$ starting from $\G^0$ and evolving in the normal direction with the velocity $\beta$. In this paper we follow the so-called direct approach in which evolution of the position vector $\vecx = \vecx(u,t)$ is governed by the equation:
\begin{equation}
\partial_t\vecx=\beta\vecN+\alpha\vecT, \quad \vecx(\cdot, 0)=\vecx^0(\cdot).
\label{eq:direct}
\end{equation}
Here $\alpha$ is the tangential component of the velocity vector.  Note that $\alpha$ has no effect on the shape of evolving closed curves, 
and the shape is determined by the value of the normal velocity $\beta$ only. 
Therefore, one can take take $\alpha\equiv0$ when analyzing analytical properties of the geometric flow driven by (\ref{eq:direct}). On the other hand, the impact of a suitable choice of a tangential velocity $\alpha$ on construction of robust and stable numerical schemes has been pointed out by many authors (see \cite{SY2008,SY2011} and references therein). 

In what follows, we shall assume that $\beta=\delta(\nu) k + \scF_\Gamma$ where $\delta(\nu)>0$ is a strictly positive $2\pi$-periodic smooth function of the tangent angle $\nu$ and $\scF_\Gamma$ is a nonlocal part of the normal velocity depending on the entire shape of the curve $\Gamma$. According to \cite{MikulaS2004a} (see also \cite{MikulaS2001, MikulaS2004b}) the system of PDEs governing evolution of plane curves evolving in the normal and tangential directions with velocities $\beta$ and $\alpha$ reads as follows:
\begin{align}
& \partial_t k=\partial_s^2 \beta+\alpha\partial_sk + k^2\beta, 
\label{eq:equation-k} \\
& \partial_t\nu= \partial_s \beta  + \alpha k,
\label{eq:equation-nu} \\
& \partial_t g=\left(-k\beta+\partial_s\alpha\right)g,
\label{eq:equation-g} \\
& \partial_t\vecx=\delta(\nu)  \partial_s^2\vecx + \alpha\partial_s\vecx +  \scF_\Gamma \vecN, 
\label{eq:equation-x}
\end{align}
for $u\in[0, 1]$ and $t>0$. Here $g=|\partial_u\vecx|$ is the so-called local length (c.f. \cite{MikulaS2001}). A solution to (\ref{eq:equation-k})--(\ref{eq:equation-x}) is subject to periodic boundary conditions for $g, k, \vecx$ at $u=0,1$, $\nu(0,t)\equiv\nu(1,t)$ mod($2\pi$) and the initial condition 
$k(\cdot,0)=k_0(\cdot), \nu(\cdot,0)=\nu_0(\cdot), g(\cdot,0)=g_0(\cdot), \vecx(\cdot,0)=\vecx_0(\cdot)$ corresponding to the initial curve $\Gamma^0=\Image(\vecx^0)$. 

Local existence and continuation of a classical smooth solution to system  (\ref{eq:equation-k})--(\ref{eq:equation-x}) has been investigated by the authors in \cite{SY2008,SY2011}. In this paper we therefore take for granted that classical solutions to (\ref{eq:equation-k})--(\ref{eq:equation-x}) exists on some maximal time interval $[0,T_{max})$ (c.f. \cite{SY2011,MikulaS2004b}).

\section{The Wulff shape and interfacial energy functional}

\begin{figure}[ht]
\begin{center}
\begin{tabular}{@{}ccc@{}}
\scalebox{0.45}{\includegraphics{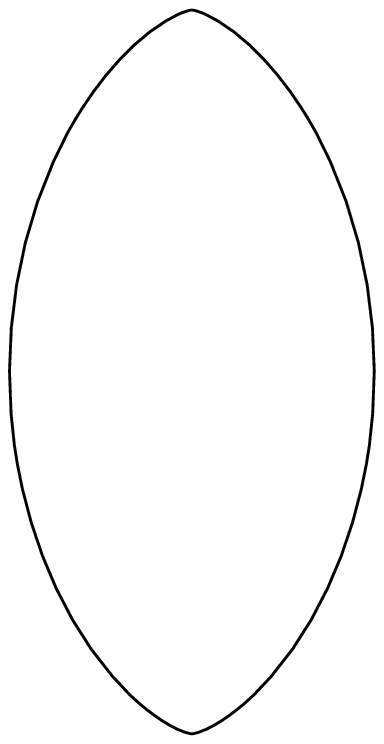}} &
\scalebox{0.45}{\includegraphics{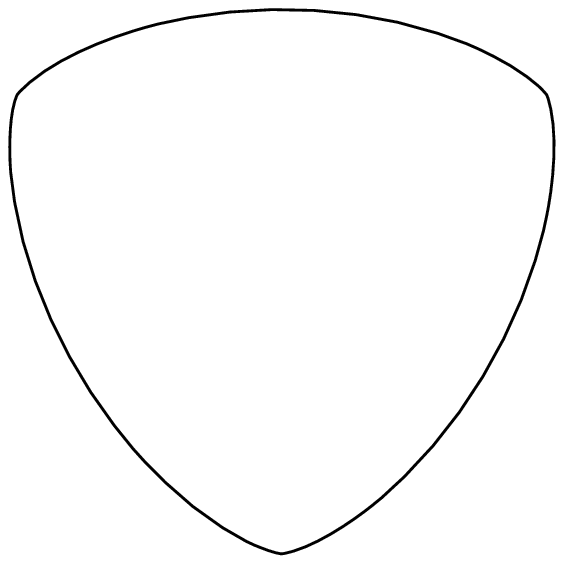}} &
\scalebox{0.45}{\includegraphics{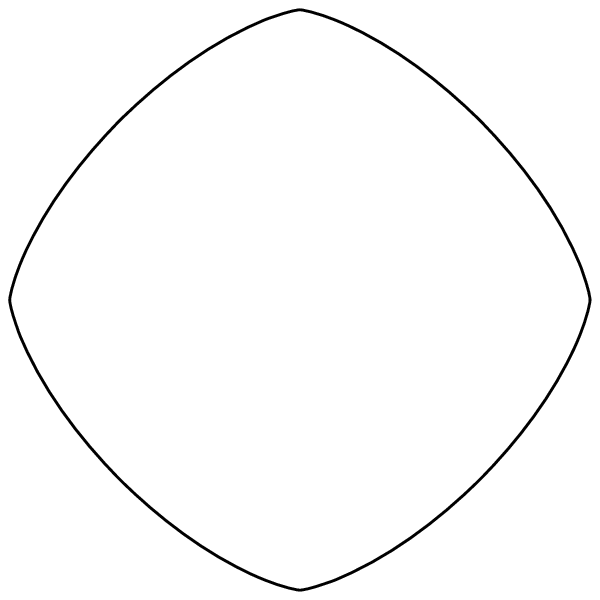}} \\
$m=2$ & $m=3$ & $m=4$ 
\end{tabular}

\begin{tabular}{@{}cc@{}}
\scalebox{0.45}{\includegraphics{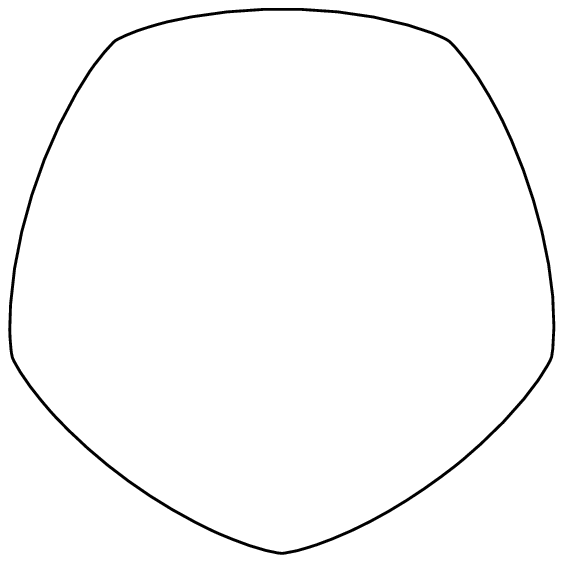}} &
\scalebox{0.45}{\includegraphics{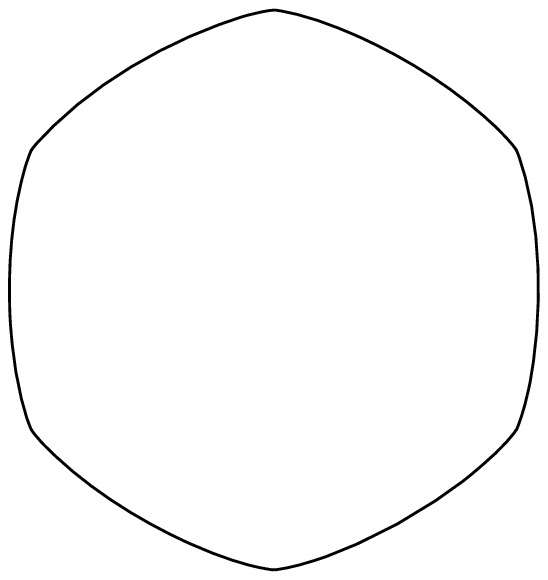}}\\
 $m=5$ & $m=6$
\end{tabular}

\end{center}
\caption{%
The Wulff shapes $W_\sigma$ for  $m=2, \cdots, 6$ and $\eps=0.99/(m^2-1)$.
}
\label{fig:Wulff}
\end{figure}

The anisotropic curvature driven flow of embedded closed plane curves is associated with the so-called interfacial energy density (anisotropy) function $\sigma$ defined on $\Gamma$. It is assumed that $\sigma=\sigma(\nu)$ is a strictly positive function depending on the tangent angle $\nu$ only. With this notation we can introduce the total interfacial energy 
\[
L_{\sigma}(\Gamma)=\int_{\Gamma}\sigma(\nu)\,\ds
\]
associated with a given anisotropy density function $\sigma$. If $\sigma\equiv1$ then  $L_1(\Gamma)$ is just the total length $L(\Gamma)$ of a curve $\Gamma$. The Wulff shape is defined as an intersection of hyperplanes: 
\[
W_{\sigma}=\bigcap_{\nu\in S^1}\biggl\{\vecx=(x_1,x_2)^{\mathrm{T}};\ -\vecx \cdot \vecN \leq \sigma(\nu)\biggr\}.
\]
If the boundary $\partial W_\sigma$ of the Wulff shape is smooth and it is parameterized by 
$\partial W_\sigma =\{ \vecx=-\sigma(\nu)\vecN+a(\nu)\vecT, \nu\in[0,2\pi] \}$, then, it follows from the relation $\partial_s\nu=k$ that 
\[
\vecT
=\partial_s\vecx
=(-\sigma'(\nu)+a(\nu))k\vecN+(\sigma(\nu)+a'(\nu))k\vecT.
\]
Hence $a(\nu)=\sigma'(\nu)$ and $(\sigma(\nu)+\sigma''(\nu))k=1$ holds and the boundary $\partial W_\sigma$ can be parameterized as follows:
\[
\partial W_{\sigma}=\left\{\vecx;\ \vecx=-\sigma(\nu)\bm{N}+\sigma'(\nu)\bm{T}, \ \nu\in [0,2\pi] \right\}, 
\]
and its curvature is given by $k=(\sigma(\nu)+\sigma''(\nu))^{-1}$. 
Let us denote by $k_\sigma$ the anisotropic curvature defined by $k_\sigma:= (\sigma(\nu)+\sigma''(\nu)) k$. It means that the anisotropic curvature $k_\sigma$ of the boundary $\partial W_\sigma$ of the Wulff shape $W_\sigma$ is constant, $k_\sigma \equiv 1$. 
Moreover, the area $|W_{\sigma}| = A(\partial W_\sigma)$ of the Wulff shape satisfies:
\begin{eqnarray*}
|W_{\sigma}|
&=&-\frac{1}{2}\int_{\partial W_{\sigma}}\vecx \cdot \vecN\,\ds
=\frac{1}{2}\int_{\partial W_{\sigma}}\sigma(\nu)\,\ds
\\
&=&\frac{1}{2}L_{\sigma}(\partial W_{\sigma}).
\end{eqnarray*}
Clearly, $|W_{1}|=\pi$ for the case $\sigma\equiv1$. If we consider the anisotropy density function $\sigma(\nu) = 1 + \varepsilon \cos(m \nu)$ for $m=2, 3, \cdots, \ \varepsilon (m^2-1)<1$ then the area of $W_\sigma$ can be easily calculated: 
\begin{eqnarray}
|W_\sigma|
&=&\frac{1}{2}\int_{\partial W_{\sigma}}\sigma(\nu)\,\ds
=\frac{1}{2}\int_0^{2\pi}\sigma(\sigma''+\sigma)\dnu
\nonumber \\
&=&\frac{\pi}{2}(2-\eps^2(m^2-1)). 
\label{eq:areaW}
\end{eqnarray}
In Fig~\ref{fig:Wulff} we plot shapes of $\partial W_\sigma$ for various degrees $m$.

Since the global quantities evaluated over the closed curve $\Gamma$ do not depend on the tangential velocity $\alpha$ 
we may take $\alpha\equiv0$. Hence $\partial_t g = -k\beta g$ and $\partial_t \nu = \partial_s\beta$. 
These identities follow from (\ref{eq:equation-nu}) and (\ref{eq:equation-g}) with $\alpha\equiv0$. 
Recall that $\partial_s \nu = k$. 
Therefore $\partial_s \sigma'(\nu) = \sigma''(\nu)\partial_s\nu = \sigma''(\nu)k$ and so $\int_\Gamma   \sigma''(\nu)k \ds = 0$. 
Hence
\[
\int_\Gamma k_\sigma \ds = \int_\Gamma \sigma k \ds
\]
holds. For the time derivative of $\int_\Gamma k_\sigma \ds$ we obtain
\begin{align*}
\frac{\mathrm{d}}{\dt} \int_\Gamma k_\sigma \ds 
&= \frac{\mathrm{d}}{\dt} \int_0^1 \sigma k g \du
= \int_0^1[ \partial_t(\sigma k) g + \sigma k \partial_t g] \du
\\
&= \int_\Gamma[ \partial_t(\sigma k)  -  \sigma k^2 \beta] \ds
\\
&= \int_\Gamma[ k \partial_t\sigma(\nu)  + \sigma(\nu)\partial_t k  -  \sigma(\nu) k^2 \beta] \ds 
\\
&= \int_\Gamma[ k \sigma'(\nu) \partial_t\nu   + \sigma(\nu)(\partial_t k  -  k^2 \beta)] \ds
\\
&= \int_\Gamma[ k \sigma'(\nu) \partial_s\beta   + \sigma(\nu) \partial_s^2 \beta] \ds= 0,
\end{align*}
because $\partial_t k  -  k^2 \beta = \partial_s^2 \beta$ and 
$k \sigma'(\nu) =  \sigma'(\nu)\partial_s\nu = \partial_s\sigma(\nu)$. 
From the previous equality we can deduce the following identity:
\begin{equation}
\int_{\Gamma^t} k_\sigma \ds = \int_{\Gamma^0} k_\sigma \ds, \quad \hbox{for any }\ 0\le t < T_{max},
\label{conservation}
\end{equation}
where the family of planar embedded closed curves $\Gamma^t, t\in [0,T_{max})$, evolves in the normal direction with the velocity $\beta$. 

Now, let us consider an evolving family of plane embedded closed curves $\Gamma^t, t\in[0,T]$, homotopicaly connecting a given curve $\Gamma=\Gamma^0$ and  the boundary $\Gamma^T=\partial W_\sigma$ of the Wulff shape $W_\sigma$. The homotopy can be realized by taking a suitable normal velocity $\beta$ (eventually depending on the position vector $\vecx$). Using such a normal velocity we deduce the identity: 
\begin{equation}
\int_\Gamma k_\sigma \ds = \int_{\partial W_\sigma} k_\sigma \ds 
= L(\partial W_\sigma).
\label{identity}
\end{equation}
It means that $\int_\Gamma k_\sigma \ds$ is equal to the length of the boundary $\partial W_\sigma$ of the Wulff shape. The same result has been recently obtained by Barrett \emph{et al.} in \cite[Lemma 2.1]{BGN}.
We can say that identity (\ref{identity}) is a generalization of the rotation number: 
$\frac{1}{2\pi}\int_\Gamma k\ds=1$, since $2\pi=L(\partial W_1)$.

\begin{remark}
Identity (\ref{identity}) can be easily shown for convex curves. Indeed, if $\Gamma$ is convex then its arc-length parameterization $s$ can be reparameterized by the tangent angle $\nu\in[0,2\pi]$. 
We have $\partial_s \nu = k >0$ and therefore $\ds = k^{-1} \dnu$. Hence
\[
\int_\Gamma k_\sigma \ds 
= \int_\Gamma \sigma k \ds 
= \int_0^{2\pi} \sigma(\nu) \dnu.
\]
For the length $L(\partial W_\sigma)$ of the boundary of a convex Wulff shape we obtain
\begin{eqnarray*}
L(\partial W_\sigma)
&=& \int_{\partial W_\sigma} \ds 
= \int_0^{2\pi} \frac{1}{k} \dnu 
\\
&=& \int_0^{2\pi} [\sigma(\nu) + \sigma''(\nu)] \dnu 
= \int_0^{2\pi} \sigma(\nu)  \dnu.
\end{eqnarray*}

Therefore $\int_\Gamma k_\sigma \ds = L(\partial W_\sigma)$ because $\int_0^{2\pi} \sigma''(\nu) \dnu = 0$ and $k=[\sigma(\nu) + \sigma''(\nu)]^{-1}$ on $\partial W_\sigma$. If $\Gamma$ is not convex we can apply the famous Grayson's theorem \cite{Gr}. We let it evolve according to the normal velocity $\beta=k$ until a time $t=T$ when $\Gamma^T$ becomes convex. Using (\ref{conservation}) and previous argument we again obtain identity (\ref{identity}). 
\end{remark}

Let us denote by $L_1$ the total interfacial energy corresponding to $\sigma\equiv 1$, i.e. $L_1\equiv L$. Let $\Gamma=\partial W_1$ be the unit circle. 
Then, by applying identity (\ref{identity}), we deduce 
\begin{equation}
L_1(\partial W_\sigma) = L_\sigma(\partial W_1). 
\label{L1Lsigma}
\end{equation}
Latter identity can be rephrased as follows: the length of the boundary $\partial W_\sigma$ of the Wulff shape equals to the total interfacial energy of the unit circle. 
It can be easily generalized to the case of arbitrary two anisotropies $\sigma(\nu)$ and $\mu(\nu)$.  We have the following proposition:

\begin{theorem}\label{prop:duality}
Let $\sigma$ and $\mu$ be two smooth anisotropy functions satisfying $\sigma(\nu)+\sigma''(\nu)>0, \mu(\nu)+\mu''(\nu)>0$. Then the duality 
\begin{equation}
L_\mu(\partial W_\sigma) = L_\sigma(\partial W_\mu)
\label{LmuLsigma}
\end{equation}
between total interfacial energies  of boundaries $\partial W_\sigma$ and $\partial W_\mu$ of Wulff shapes holds. 
\end{theorem}

\noindent{P r o o f.}
Notice that the Wulff shapes $W_\sigma$ and $W_\mu$ are convex sets because 
$\sigma(\nu)+\sigma''(\nu)>0$ and $\mu(\nu)+\mu''(\nu)>0$ hold. 
For the curvature $k$ at the boundary $\partial W_\sigma$ we have $k =[\sigma(\nu)+\sigma''(\nu)]^{-1}$ and so
\begin{align}
\label{LmuLsigmaProof}
L_\mu(\partial W_\sigma) 
&= \int_{\partial W_\sigma} \mu(\nu) \ds 
\\
&= \int_0^{2\pi} \mu(\nu) \frac{1}{k} \dnu 
= \int_0^{2\pi} \mu(\nu) (\sigma(\nu)+\sigma''(\nu)) \dnu  \nonumber
\\
&= \int_0^{2\pi}[ \mu(\nu) \sigma(\nu)- \sigma'(\nu) \mu'(\nu)] \dnu
= L_\sigma(\partial W_\mu), 
\end{align}
arguing vice versa. \hfill $\diamondsuit$

\section{Gradient flow for the anisoperimetric ratio.}
\label{sec:nonlocal}

Recall that for the enclosed area $A=A(\Gamma^t)$ and the total length $L=L(\Gamma^t)$ for a flow of embedded closed plane curves driven in normal direction by the velocity $\beta$ we have
\begin{equation}
\frac{\mathrm{d}}{\dt} A + \int_{\Gamma^t} \beta \ds = 0,\quad
\frac{\mathrm{d}}{\dt} L + \int_{\Gamma^t} k \beta \ds = 0,
\label{eq:area}
\end{equation}
(c.f. \cite{MikulaS2001}). Using governing equations (\ref{eq:equation-k})--(\ref{eq:equation-x}), for the total interfacial energy $L_\sigma = L_\sigma(\Gamma^t)$ of a curve $\Gamma^t$,  we obtain 
\begin{align}
\label{eq:lengthaniso}
\frac{\mathrm{d}}{\dt} L_{\sigma} 
&= \frac{\mathrm{d}}{\dt} \int_{\Gamma}\sigma(\nu)\ds
= \frac{\mathrm{d}}{\dt} \int_0^1\sigma(\nu) g \du
\\
&= \int_0^1 [ \sigma'(\nu)\partial_t\nu  g + \sigma(\nu) \partial_tg ] \du
\nonumber
\\
&= \int_{\Gamma} [ \sigma'(\nu)\partial_s \beta  - \sigma(\nu) k\beta ] \ds
\\
&= -\int_{\Gamma} [ \sigma''(\nu)\partial_s \nu \beta  +  \sigma(\nu) k\beta ] \ds 
\\
&= - \int_{\Gamma} [\sigma''(\nu) +  \sigma(\nu) ] k\beta \ds 
= - \int_{\Gamma}  k_\sigma \beta  \ds.
\nonumber
\end{align}
Here we have used the governing equations (\ref{eq:equation-g}) and (\ref{eq:equation-nu}) (with $\alpha\equiv 0$) and the identity $\partial_s\nu = k$. For the anisoperimetric ratio 
\[
\Pi_{\sigma}(\Gamma)=\frac{L_{\sigma}(\Gamma)^2}{4|W_{\sigma}|A(\Gamma)},
\]
we have $\Pi_{\sigma}(\Gamma) \ge 1$ and, in particular, $\Pi_{\sigma}(\partial W_{\sigma})=1$ (see Remark~\ref{remaniso}). 
Taking into account identities (\ref{eq:lengthaniso}) and (\ref{eq:area}) we obtain
\begin{eqnarray*}
\frac{\mathrm{d}}{\dt}\Pi_\sigma 
&=& \frac{L_\sigma \partial_t L_\sigma }{2|W_{\sigma}| A} - \frac{L_\sigma^2 \partial_t A}{4 |W_{\sigma}| A^2}
\\
&=& - \frac{L_\sigma}{2 |W_{\sigma}| A} \int_{\Gamma}\left(k_\sigma  - \frac{L_\sigma }{2 A}\right) \beta \ds.
\end{eqnarray*}

Hence, the flow driven in the normal direction by the non-locally dependent velocity
\begin{equation}
\beta = k_\sigma - \frac{L_\sigma}{2A} 
\label{eq:isoperimnormalaniso}
\end{equation}
represents a gradient flow for the anisoperimetric ratio $\Pi_\sigma$ with the property  $\partial_t \Pi_\sigma  <0$ for $\beta\not\equiv0$. Notice that $\beta\equiv0$ on $\Gamma$ if and only if $\Gamma \propto \partial W_\sigma$, i.e. $\Gamma$ is homotheticaly similar to $\partial W_\sigma$.

In the case $\sigma\equiv 1$ the isoperimetric ratio gradient flow has been analyzed by Jiang and Zhu in \cite{Jiang2008} and by the authors in \cite{SY2008}. In this case the normal velocity has the form: 
$\beta = k - L/(2 A)$.

\section{A mixed anisoperimetric inequality}

The aim of this section is to prove a mixed anisoperimetric inequality of the form
\begin{equation}
\frac{L_\sigma(\Gamma) L_\mu(\Gamma)}{A(\Gamma)} \ge K_{\sigma,\mu},
\label{mixedanisoperim}
\end{equation}
which holds for any $C^2$ smooth Jordan curve $\Gamma$ in the plane. 
Here $K_{\sigma,\mu}>0$ is a constant depending only on the anisotropy functions $\sigma$ and $\mu$ such that 
$\sigma(\nu)+\sigma''(\nu)>0$ and $\mu(\nu)+\mu''(\nu)>0$ hold for any $\nu$. The existence of a minimizer of the mixed anisoperimetric ratio is discussed in Remark~\ref{minimizer}. The idea of the proof of the inequality (\ref{mixedanisoperim}) is rather simple and consists in solving the constrained minimization problem:
\begin{equation}
\min_{\Gamma} L_\sigma(\Gamma), \qquad \hbox{s.t.}\ \ L_\mu(\Gamma) = c A(\Gamma),
\label{minimproblem}
\end{equation}
where $c>0$ is a given constant. To this end, let us assume that a curve $\Gamma=\Gamma(\vecx)$ is parameterized by a $C^2$ smooth function $\vecx: S^1\to \R^2$. If we denote $g\equiv g(\vecx) = |\partial_u \vecx|$ the local length then, for the derivative of $g$ in the direction $\vecy:S^1\to\R^2$, 
we obtain  $2 g(\vecx) g'(\vecx) \vecy = 2 (\partial_u \vecx \cdot \partial_u \vecy)$ and so
\begin{equation}
g'(\vecx) \vecy = (\vecT \cdot \partial_s \vecy) g.
\label{derg} 
\end{equation}
Here and here after, for scalar-valued function $f(\vecx)$ and vector-valued function $\vecf(\vecx)=(f_1(\vecx), f_2(\vecx))^{\mathrm{T}}$ 
we denote their derivatives in the direction $\vecy$ by
\[
f'(\vecx)\vecy:=\nabla f(\vecx)\cdot\vecy=\lim_{\eps\to 0}\frac{f(\vecx+\eps\vecy)-f(\vecx)}{\eps}, 
\]
\[
\vecf'(\vecx)\vecy:=\begin{pmatrix}f_1'(\vecx)\vecy\\ f_2'(\vecx)\vecy\end{pmatrix},  
\]
respectively. 

As for the tangent vector $\vecT=\vecT(\vecx)= (\cos\nu, \sin\nu)^{\mathrm{T}}$ we have $\vecT(\vecx) = g^{-1}\partial_u\vecx$ and so 
$\vecT'(\vecx) \vecy 
= g^{-1} \partial_u \vecy  - g^{-2} \partial_u\vecx \, g'(\vecx) \vecy 
= \partial_s\vecy -(\vecT \cdot \partial_s \vecy) \vecT 
= (\vecN \cdot \partial_s \vecy) \vecN$. 
As $\vecN= (-\sin\nu, \cos\nu)^{\mathrm{T}}$, for the derivative of the tangent angle $\nu=\nu(\vecx)$, we obtain  
\begin{equation}
\nu'(\vecx) \vecy = \vecN \cdot \partial_s \vecy.
\label{dernu} 
\end{equation}
Recall that $k_\sigma :=(\sigma(\nu)+ \sigma''(\nu)) k$ and $\partial_s\nu = k$.  Since $L_\sigma(\Gamma)= \int_\Gamma \sigma \ds = \int_0^1 \sigma(\nu) g \du$ 
we obtain
\begin{align*}
L_\sigma^\prime (\Gamma(\vecx)) \vecy 
& = \int_0^1 
\left[
\sigma'(\nu) \nu'(\vecx) \vecy g + \sigma (\nu) g'(\vecx) \vecy
\right] \du 
\\
& = \int_\Gamma 
\left[ 
\sigma'(\nu) (\vecN \cdot \partial_s\vecy) + \sigma(\nu) (\vecT \cdot \partial_s\vecy) 
\right]\ds
\\
& = - \int_\Gamma 
\left[
\sigma''(\nu)\partial_s\nu  (\vecN \cdot \vecy) 
- \sigma'(\nu) k (\vecT \cdot \vecy)\right.
\\
& \left. \hskip 1truecm + \sigma'(\nu) \partial_s\nu  (\vecT \cdot \vecy)
+ \sigma(\nu) k (\vecN \cdot \vecy) 
\right]\ds
\\
& 
= - \int_\Gamma 
\left[ 
\sigma(\nu)+ \sigma''(\nu) \right] k (\vecN \cdot \vecy) \ds 
\\
&= - \int_\Gamma k_\sigma  (\vecN \cdot \vecy)  \ds .
\end{align*}
Hence
\begin{eqnarray}
L_\sigma^\prime(\Gamma(\vecx)) \vecy 
&=& - \int_\Gamma  k_\sigma (\vecN \cdot \vecy) \ds, \quad 
\nonumber \\
L_\mu^\prime(\Gamma(\vecx)) \vecy  &=& - \int_\Gamma  k_\mu (\vecN \cdot \vecy) \ds.
\label{derL}
\end{eqnarray}
For the area $A=A(\Gamma)$  enclosed by a Jordan curve $\Gamma=\Gamma(\vecx)$ we have
$A =  \frac12 \int_0^1 \det(\vecx, \partial_u\vecx ) \du$. Therefore
\begin{align*}
A'(\Gamma(\vecx))\vecy &= \frac12 \int_0^1 \det( \vecy, \partial_u\vecx ) + \det(\vecx, \partial_u\vecy ) \du 
\\
&=  \int_0^1 \det(\vecy, \partial_u\vecx ) \du= \int_\Gamma \det(\vecy, \vecT)\ds.
\end{align*}

Since $\det(\vecy, \vecT) = - \vecy \cdot \vecN$ we obtain 
\begin{equation}
A'(\Gamma(\vecx)) \vecy = - \int_\Gamma \vecN \cdot \vecy \ds.
\label{derA} 
\end{equation}
In order to solve the constrained minimization problem (\ref{minimproblem}) we introduce the  Lagrange function 
${\mathcal L}(\vecx,\lambda) = L_\sigma (\Gamma(\vecx)) + \lambda (L_\mu (\Gamma(\vecx)) - c A(\Gamma(\vecx)))$ with $\lambda>0$. 

Then the first order condition for $\bar\Gamma =\Gamma(\bar\vecx)$ to be a minimizer of  (\ref{minimproblem}) reads as follows: 
$ 0= {\mathcal L}^\prime_{\vecx}(\vecx,\lambda) \vecy \equiv  L_\sigma^\prime(\Gamma(\vecx)) \vecy + \lambda (L_\mu^\prime(\Gamma(\vecx)) \vecy -c  A'(\Gamma(\vecx)) \vecy)$ at $\vecx=\bar\vecx$. Latter equality has to be satisfied for any smooth function $\vecy:S^1\to \R^2$. Taking into account (\ref{derL}) and (\ref{derA}) we obtain 
\[
k_\sigma + \lambda k_\mu = \lambda c, \quad\hbox{on}\ \ \bar\Gamma. 
\]
It means that 
\begin{equation}
k_{\bar\sigma} = \lambda c, \quad\hbox{on}\ \ \bar\Gamma, \quad \hbox{where}\ \ 
\bar\sigma = \sigma + \lambda\mu.
\label{neccondk}
\end{equation}
In other words, $\bar\Gamma = \frac{1}{\lambda c} \partial W_{\bar\sigma}$ (up to an affine translation in the plane $\R^2$). 
The Lagrange multiplier $\lambda\in\R$ can be computed from the constraint $L_\mu(\bar\Gamma) = c A(\bar\Gamma)$. 
It follows from duality (\ref{LmuLsigma}) (see Proposition~\ref{prop:duality}) that 
\begin{eqnarray*}
L_\mu(\partial W_{\bar\sigma}) 
&=&  L_{\bar\sigma}(\partial W_\mu) 
= L_\sigma (\partial W_\mu) + \lambda L_\mu (\partial W_\mu)
\\
&=& L_\sigma (\partial W_\mu) + 2 \lambda A(\partial W_\mu).
\end{eqnarray*}
To calculate the enclosed area $A(\bar\Gamma) = \frac{1}{\lambda^2 c^2} A(\partial W_{\bar\sigma})$ we make use of the identity $A(\partial W_{\bar\sigma}) = \frac12 L_{\bar\sigma}(\partial W_{\bar\sigma})$. Clearly, as $\bar\sigma=\sigma+\lambda\mu$ we obtain 
\begin{align*}
L_{\bar\sigma}(\partial W_{\bar\sigma}) 
& = L_\sigma(\partial W_{\bar\sigma}) + \lambda L_\mu(\partial W_{\bar\sigma}) 
\\
& = L_{\bar\sigma}(\partial W_\sigma) + \lambda L_{\bar\sigma}(\partial W_\mu)
\\
& 
= L_\sigma(\partial W_\sigma) + \lambda L_\mu(\partial W_\sigma)
+ \lambda L_\sigma(\partial W_\mu) 
\\ 
& + \lambda^2 L_\mu(\partial W_\mu)
\\
& = 2 A(\partial W_\sigma) + 2 \lambda L_\sigma(\partial W_\mu) + 2 \lambda^2 A(\partial W_\mu).
\end{align*}
Since $\frac{1}{\lambda c} L_\mu(\partial W_{\bar\sigma}) = L_\mu(\bar\Gamma) = c A(\bar\Gamma) = \frac{c}{\lambda^2 c^2} A(\partial W_{\bar\sigma})$ we end up with the identity
\begin{eqnarray*}
&&\frac{1}{\lambda c} 
\left(
L_\sigma(\partial W_\mu) + 2\lambda A(\partial W_\mu)
\right)
\\
&=& 
\frac{c}{\lambda^2 c^2} 
\left(
 A(\partial W_\sigma) +  \lambda L_\sigma(\partial W_\mu) +  \lambda^2 A(\partial W_\mu)
\right). 
\end{eqnarray*}

Since the Lagrange multiplier $\lambda>0$ it is given by $\lambda = \sqrt{A(\partial W_\sigma)/A(\partial W_\mu)}$. Furthermore, 
\begin{eqnarray*}
L_\sigma(\partial W_{\bar\sigma}) 
&=&  L_{\bar\sigma}(\partial W_\sigma) 
= L_\sigma (\partial W_\sigma) + \lambda L_\mu (\partial W_\sigma)
\\
&=& 2  A(\partial W_\sigma) + \lambda L_\sigma (\partial W_\mu).
\end{eqnarray*}

Now, let $\Gamma$ be an arbitrary $C^2$ smooth Jordan curve in the plane. Set $c= L_\mu(\Gamma)/A(\Gamma)$. Then
\begin{align*}
\frac{L_\sigma(\Gamma) L_\mu(\Gamma)}{A(\Gamma)} 
&= c L_\sigma(\Gamma) \ge  c L_\sigma(\bar\Gamma) =\frac{c}{\lambda c} L_\sigma(\partial W_{\bar\sigma})
\\
&= 2 \sqrt{A(\partial W_\sigma) A(\partial W_\mu) } + L_\sigma (\partial W_\mu).
\end{align*}

\begin{remark}\label{minimizer}
The proof of existence of a minimizer of the mixed anisoperimetric ratio $L_\sigma(\Gamma) L_\mu(\Gamma)/A(\Gamma)$ is as follows: let $\Gamma^n=\Gamma(\vecx^n)$ be a sequence of Jordan curves minimizing this ratio. As $L_\sigma(\gamma\Gamma)=\gamma L_\sigma(\Gamma)$, and $A(\gamma\Gamma)=\gamma^2A(\Gamma)$ for each $\gamma>0$, without lost of generality, we may assume $L(\Gamma^n)=1$ for all $n\in\N$. We can also fix the barycenter of $\Gamma^n$ at the origin. 
Since $c_0 L(\Gamma)\le L_\sigma(\Gamma)\le c_1 L(\Gamma)$ where $0<c_0=\min_\Gamma\sigma\le c_1=\max_\Gamma\sigma<\infty$, then, by the isoperimetric inequality, the value of the infimum is positive. Moreover, the parameterization $\vecx^n$ of $\Gamma^n$ can be chosen in such a way that $|\partial_u \vecx^n|=L(\Gamma^n)=1$. As a consequence, the position vectors $\{\vecx^n(u), u\in[0,1]\}$ are uniformly bounded. By the Arzel\`a-Ascoli theorem there is a convergent subsequence converging to some function $\{\vecx(u), u\in[0,1]\}$ which is the minimizer of the mixed anisoperimetric ratio.

\end{remark}

In summary, we have shown the following mixed anisoperimetric inequality:

\begin{theorem}\label{anisoineq}
Let $\Gamma$ be a $C^2$ smooth Jordan curve in the plane. Then
\begin{equation}
\frac{L_\sigma(\Gamma) L_\mu(\Gamma)}{A(\Gamma)} \ge K_{\sigma,\mu},
\label{mixedanisoperim2}
\end{equation}
where $K_{\sigma,\mu} = 2 \sqrt{|W_\sigma| |W_\mu| } + L_\sigma (\partial W_\mu)$.The equality in (\ref{mixedanisoperim2}) holds if and only if the curve $\Gamma$ is homothetically similar to the boundary $\partial W_{\widetilde{\sigma}}$ of a Wulff shape corresponding to the mixed anisotropy function $\widetilde{\sigma} = \sqrt{|W_\mu|}\, \sigma + \sqrt{|W_\sigma|} \,\mu$.

\end{theorem}

\medskip

\begin{remark} \label{remaniso}
If $\sigma=\mu\equiv1$ we obtain the well known isoperimetric inequality $L(\Gamma)^2/A(\Gamma) \ge K_{1,1} \equiv 2 \sqrt{\pi^2} + L(W_1) = 4\pi$. 
If $\sigma=\mu$ we obtain the anisoperimetric inequality 
$L_\sigma(\Gamma)^2/A(\Gamma) \ge K_{\sigma,\sigma} = 2 \sqrt{|W_\sigma|^2} + L_\sigma(\partial W_\sigma) = 4 |W_\sigma|$. 
Finally, if $\mu\equiv 1$ we obtain the mixed anisoperimetric inequality
\[
\frac{L_\sigma(\Gamma) L(\Gamma)}{A(\Gamma)} \ge K_{\sigma,1} \equiv 2 \sqrt{\pi |W_\sigma|} 
+ L(\partial W_\sigma).
\]
\end{remark}

\begin{remark} \label{remaniso2}
In the case $\mu=\sigma$, the anisoperimetric inequality in the plane has been stated in a paper by G. Wulff \cite{Wulff1901} from 1901. Later, it was proved by Dinghas in \cite{Dinghas1944} for a special class of polytopes. Recently, Fonseca and M\"uller \cite{Fonseca1991}  proved the anisotropic inequality in the plane. Later Fusco \emph{et al.} \cite{Fusco2005} proved it in arbitrary dimension. Giga in \cite{Giga2003} pointed out that the anisotropic inequality where $\mu=\sigma$ are $\pi$-periodic function is the isoperimetric inequality in a suitable Minkowski metric. It is a useful tool in the proof of anisotropic version of the so-called Gage's inequality (c.f. \cite[Corollary 4.3]{Gage1993}). 

However, in all aforementioned proofs, the surface energy was associated with a functional $L_\Phi(\Gamma)= \int_\Gamma \Phi(\vecN) \ds$ where $\Phi:\R^2 \to\R$ is an absolute homogeneous anisotropy function of degree one, i.e. $\Phi(t\vecx) = |t|\Phi(\vecx)$ for any $t\in\R, \vecx\in\R^2$. The relation between our description of anisotropy and the latter one is: 
$\sigma(\nu) = \Phi(-\sin\nu, \cos\nu)$ and, conversely, $\Phi(\vecx)= \sigma(\nu)$ where $\vecx/\Vert\vecx\Vert = (-\sin\nu, \cos\nu)$. Since we do not require $\pi$-periodicity of $\sigma$, in our approach of description of anisotropy we therefore allow for non-symmetric anisotropies, like e.g. functions $\sigma$ with odd degree $m$ (see Fig \ref{fig:Wulff}) corresponding thus to anisotropy function $\Phi$ which are positive homogeneous only, i.e. $\Phi(t\vecx) = t\Phi(\vecx)$ for any $t\ge 0, \vecx\in\R^2$.

In the case of general anisotropy functions $\mu\not\equiv\sigma$, the mixed anisoperimetric inequality derived in Theorem~\ref{anisoineq} is, to our best knowledge, new even in the case of symmetric ($\pi$-periodic) anisotropy functions.  

\end{remark}

\section{Convexity preservation. Temporal area and length behavior}

In this section we analyze behavior of the enclosed area $A(\Gamma^t)$ of a curve $\Gamma^t$ evolved in the normal direction by the  anisoperimetric ratio gradient flow, i.e. $\beta = k_\sigma - L_\sigma/(2 A)$. 

First we prove the preservation of convexity result stating that the anisoperimetric ratio gradient flow preserves convexity of evolved curves. 
In the case of the isoperimetric ratio gradient flow of convex curves with $\beta = k - L/(2 A)$, the convexity preservation has been shown by Jiang and Pan in \cite{Jiang2008}. However, similarly as Mu and Zhu in \cite{Ma2012}, they utilized the Gauss parameterization of the curvature equation (\ref{eq:equation-k}) by the tangent angle $\nu$ and this is why their results are applicable to evolution of convex curves only. In our paper we first prove convexity preservation based on the analysis of the curvature equation (\ref{eq:equation-k}) with arc-length parameterization. Moreover, we show the anisoperimetric ratio gradient flow may initially increase the total length and decrease the enclosed area. This phenomenon cannot be found in the isoperimetric ratio gradient flow (c.f. \cite{Jiang2008,SY2008}). 

\begin{theorem}\label{th:convexity} Let $\Gamma^t, t\in[0,T_{max})$, be the anisoperimetric ratio gradient flow of smooth Jordan curves in the plane evolving in the normal direction by the velocity $\beta=k_\sigma - \frac{L_\sigma}{2 A}$. If the curve $\Gamma^{t_0}$ is convex at some time $t_0\in[0,T_{max})$ then $\Gamma^{t}$ remains convex for any $t\in [t_0,T_{max})$.
\end{theorem}

\noindent{P r o o f.}
Since $\partial_t k = \partial_s^2\beta + k^2 \beta$, $\partial_t\nu = \partial_s\beta = \partial_s k_\sigma$ and $\beta = k_\sigma +\scF_\Gamma$ we have
\begin{align*}
\partial_t k_\sigma 
&= \delta(\nu) \partial_t k + \delta^\prime(\nu) k \partial_t \nu 
\\
&= \delta(\nu) \partial_s^2 k_\sigma + \delta(\nu) k^2 \beta +  \delta^\prime(\nu) k \partial_s k_\sigma \\
&= \delta(\nu) \partial_s^2 k_\sigma + \frac{1}{\delta(\nu)} k_\sigma^2 \beta + \delta^\prime(\nu) k \partial_s k_\sigma,
\end{align*}
where $\delta(\nu) := \sigma(\nu) + \sigma^{\prime\prime}(\nu)>0$. 

Let us denote by $K(t) = \min_{\Gamma^t} k_\sigma(.,t)$ the minimum of the anisotropic curvature $k_\sigma = \delta(\nu) k$ over the curve $\Gamma^t$. Denote by $s^*(t)\in [0, L^t]$ the argument of the minimum of $k_\sigma$, i.e. $K(t) = k_\sigma(s^*(t), t)$. Then $\partial_s k_\sigma (s^*(t), t)=0$ and 
$\partial_s^2 k_\sigma (s^*(t), t)\ge 0$. Hence
\[
K'(t) \ge \frac{1}{\Delta(t)}K(t)^2 (K(t) +  \scF_\Gamma^t),
\]
where $\Delta(t) = \delta(\nu(s^*(t), t))$, and $\scF_\Gamma^t=-L^t_\sigma/(2A^t)$. 
Notice that $\Delta(t)\ge \Delta_{min}:=\min_{\nu} \delta(\nu)>0$ for all $t\in[0,T_{max})$. Suppose that $K$ is a solution to this  ordinary differential inequality existing on some interval $[t_0, T_{max})$ and such that $K(t_0)>0$. Then, it should be obvious that $K(t)>0$ for $t\in[t_0, T_{max})$ provided that 
\begin{equation}
\inf_{t_0\le t \le t^*}\scF_\Gamma^t > -\infty,
\label{Fbound}
\end{equation}
for every $0<t^*<T_{max}$.  In order to prove convexity preservation for the anisoperimetric ratio  gradient flow it is therefore sufficient to verify that the nonlocal part 
$\scF_\Gamma^t=-L^t_\sigma/(2A^t)$ remains bounded from below for $t\ge t_0$. To prove boundedness of $\scF_\Gamma^t$ from below  we utilize a property of the anisoperimetric ratio. Indeed, as $\beta = k_\sigma -L_\sigma/(2A)$ represents gradient flow for the anisoperimetric ratio 
$\Pi^t_\sigma = (L_\sigma^t)^2/(4|W_\sigma| A^t)$ we have $1\le \Pi^t_\sigma \le \Pi^0_\sigma$ for all $t\in [0,T_{max})$. Thus 
\[
\scF_\Gamma^t=-\frac{L^t_\sigma}{2 A^t} \ge - \frac{(L^0_\sigma)^2}{2 A^0} \frac{1}{L^t_\sigma}.
\]
Now, since the classical solution exists on the time interval $[0,T_{max})$ then $\inf_{0\le t\le t^*} L^t_\sigma >0$ for each $0<t^*<T_{max}$ and the estimate (\ref{Fbound}) follows.  \hfill $\diamondsuit$

\bigskip 

In what follows, we shall investigate the enclosed area and length behavior of curves evolved by the normal velocity $\beta =k_\sigma - L_\sigma/(2A)$ representing thus a gradient flow for the anisoperimetric ratio. 

Using the area equation (\ref{eq:area}) we obtain 
\begin{eqnarray}
\frac{\mathrm{d}}{\dt} A 
&=& - \int_{\Gamma} \beta \ds 
= - \int_{\Gamma}\left(k_\sigma - \frac{L_\sigma}{2 A}\right)\ds
\nonumber
\\
&=& - L (\partial W_\sigma) + \frac{L L_\sigma}{2 A}.
\label{der1A}
\end{eqnarray}
By applying the isoperimetric inequality (see Remark~\ref{remaniso}) for the case $\sigma\equiv 1$ and any curve $\Gamma=\Gamma^t$, the following inequality:
\[
\frac{\mathrm{d}}{\dt} A = - 2\pi  + \frac{L(\Gamma)^2}{2 A(\Gamma)} \ge 0
\]
holds. It means that the gradient flow for the isoperimetric ratio does not decrease the enclosed area. On the other hand, if anisotropy density function $\sigma\not\equiv const$,  then for a curve $\Gamma=\Gamma^t \propto \partial W_{\bar\sigma}$ corresponding to the Wulff shape 
$W_{\bar\sigma}$ with the anisotropy function $\bar\sigma = \sqrt{\pi}\, \sigma + \sqrt{|W_\sigma|}$ we obtain 
\begin{eqnarray}
\frac{\mathrm{d}}{\dt} A 
&=& - L (\partial W_\sigma) + \frac{L L_\sigma }{2 A}
= - L (\partial W_\sigma) + \frac{1}{2} K_{\sigma,1} 
\nonumber 
\\
&=& \sqrt{\pi |W_\sigma|} 
- \frac12 L(\partial W_\sigma) < 0
\label{derAneg}
\end{eqnarray}
due to the isoperimetric inequality $L(\partial W_\sigma)^2 \ge 4 \pi A(\partial W_\sigma)=4 \pi |W_\sigma|$. 
It means that the gradient flow for the anisoperimetric ratio may initially decrease the enclosed area for special initial curves. 

Next we recall the isoperimetric inequality by Gage. According to \cite{Gage1983} the following inequality holds:
\begin{equation}
\int_\Gamma k^2 \ds 
\ge 
\pi \frac{L}{A}
\label{Gageaniso}
\end{equation}
for any  convex $C^2$ smooth Jordan curve in the plane. The equality in (\ref{Gageaniso}) holds iff $\Gamma$ is a circle. Therefore, in the case of isoperimetric gradient flow with $\sigma\equiv 1$ and the convex curve $\Gamma^t$, we have 
\begin{equation}
\frac{\mathrm{d}}{\dt} L^t=  - \int k\beta \ds  = - \int_\Gamma k^2 \ds + \pi\frac{L}{A} \le 0.
\label{length-decrease}
\end{equation}
However, if $\sigma\not\equiv const$ is a smooth nonconstant anisotropy density function such that $\sigma +\sigma^{\prime\prime}>0$, there exists an initial curve $\Gamma^0$ such that the length $L^t$ may initially increase, i.e. $\frac{\mathrm{d}}{\dt} L^t>0$ at $t=0$. Indeed, let $\Gamma^0$ be an initial curve which homothetically similar to the boundary $\partial W_{\bar\sigma}$ of the Wulff shape corresponding to the mixed anisotropy function $\bar\sigma = a \sigma + b$ where $a, b>0$ are constants. Then
\[
k_\sigma = k_{\frac{\bar\sigma -b}{a}} = \frac{1}{a} k_{\bar\sigma}   - \frac{b}{a} k.
\]
Hence $k_\sigma = \frac{1}{a} - \frac{b}{a} k $ on the Wulff shape $\Gamma=\Gamma^0 = \partial W_{\bar\sigma}$ because $k_{\bar\sigma} \equiv 1$ on $\partial W_{\bar\sigma}$. 
Using (\ref{eq:area}), we have
\begin{eqnarray*}
\frac{\mathrm{d}}{\dt} L^t
&=&  - \int_\Gamma k k_\sigma  \ds   +  \pi \frac{L_\sigma(\partial W_{\bar\sigma})}{A(\partial W_{\bar\sigma})}  
\\
&=& 
\frac{b}{a}  \int_\Gamma k^2 \ds - \frac{2\pi}{a} + \pi\frac{L_\sigma(\partial W_{\bar\sigma} )}{A(\partial W_{\bar\sigma})}
\end{eqnarray*}
at $t=0$. Since $\sigma=(\bar\sigma -b)/a$ we have
\begin{align*}
L_\sigma(\partial W_{\bar\sigma}) &= \int_{\partial W_{\bar\sigma}} \frac{\bar\sigma -b}{a} \ds
\\
&= \frac{1}{a} L_{\bar\sigma}(\partial W_{\bar\sigma})  - \frac{b}{a} L(\partial W_{\bar\sigma})
\\
&= \frac{2}{a} A(\partial W_{\bar\sigma})   - \frac{b}{a} L(\partial W_{\bar\sigma}).
\end{align*}

Thus 
\[
\frac{\mathrm{d}}{\dt} L^t= \frac{b}{a} \left( 
 \int_{\partial W_{\bar\sigma}} k^2 \ds -  \pi\frac{L(\partial W_{\bar\sigma})}{A(\partial W_{\bar\sigma})} \right)  >0 \quad\hbox{at}\ t=0, 
\]
due to inequality (\ref{Gageaniso}) and the fact that $\partial W_{\bar\sigma}$ is a convex curve different from a circle for $\bar\sigma\not\equiv const$. 

\medskip

In summary, we have shown the following result.

\begin{theorem}\label{th:areadecrasing} 
If $\sigma\equiv 1$ then the isoperimetric ratio gradient flow with the normal velocity $\beta=k-L/(2A)$ is area nondecreasing and length nonincreasing flow  of smooth Jordan curves $\Gamma^t, t\in[0,T_{max})$ in the plane provided that $\Gamma^0$ is a convex curve. 

Assume the anisotropy function $\sigma$ is not constant and such that $\sigma + \sigma^{\prime\prime}>0$. Let $\Gamma^{0}$ be an initial curve which is homothetically similar to the boundary $\partial W_{\bar\sigma}$ of a Wulff shape with the modified anisotropy density function $\bar\sigma = a \sigma + b$ where $a,b$ are constants, $a,b>0$. Then, for the anisoperimetric ratio gradient flow  $\Gamma^t, t\in[0,T_{max})$, evolving in the normal direction by the velocity $\beta=k_\sigma - \frac{L_\sigma}{2 A}$, we have

\begin{enumerate}
\item $\frac{\mathrm{d}}{\dt} L(\Gamma^t) > 0$ at $t=0$. 

\item If, moreover, $a/b = \sqrt{\pi}/\sqrt{|W_\sigma|}$ then $\frac{\mathrm{d}}{\dt} A(\Gamma^t) < 0$ at $t=0$.

\end{enumerate}

\end{theorem}

\section{A counterexample to the comparison principle}

The aim of this section is to demonstrate that the comparison property does not hold under the anisoperimetric gradient flow, which is quite in a contrast to the total-length gradient flow $\beta=k$. It is a well-known fact that the comparison argument plays a key role in the proof of the famous Gage-Hamilton-Grayson theorem
for the curvature driven flow $\beta=k$, and it states that two smooth curves, one of them included in the closure of the interior of the second one, evolved by the normal velocity $\beta=k$ never intersects each other \cite{GH1986,Gr}. The aim of this section is to show that the analogous comparison property does not hold for the anisoperimetric gradient flow. As the flow $\beta=k_\sigma - L_\sigma/(2A)$ is nonlocal, violation of comparison principle can be expected. Nevertheless, we provide an explicit construction of a counterexample in this section. 

Clearly, any curve $\Gamma$ homotheticaly similar to the Wulff shape $\partial W_\sigma$ is a stationary curve, i.e. $\beta\equiv 0$ on $\Gamma$. Indeed, $k_\sigma \equiv 1$ on $\partial W_\sigma$ and $\scF_{\partial W_\sigma} = - L_\sigma(\partial W_\sigma)/(2 A(\partial W_\sigma)) = -1$ and therefore $\beta=k_\sigma + \scF_{\partial W_\sigma} \equiv 0$ on $\partial W_\sigma$.

\begin{figure}[ht]
\centering
\includegraphics[width=1.5in]{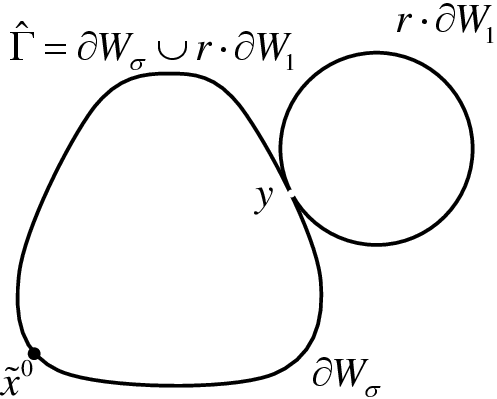}
\includegraphics[width=1.5in]{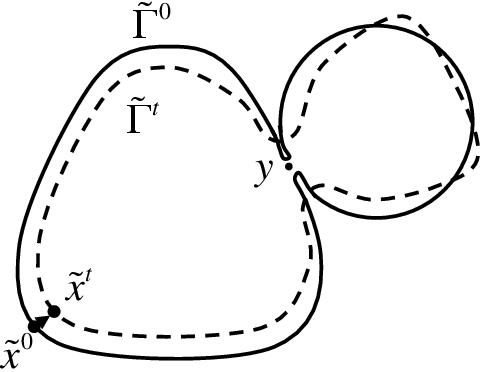} 
\caption{%
An initial nonsmooth curve $\hat\Gamma$ (left) and its smooth perturbation $\tilde\Gamma^0$ (right). Failure of a comparison principle occurs at the point $\tilde\vecx^0$ and $0<t\ll 1$. 
}
\label{fig:failurecomparison}
\end{figure}

In what follows, we shall construct a smooth initial curve $\tilde\Gamma^0$ containing in its interior the Wulff shape $\partial W_\sigma$ and such that $\tilde\Gamma^t$ intersects the stationary Wulff shape $\partial W_\sigma$ for all sufficiently small times $0<t\ll 1$. The construction is as follows. First, we shall  construct a nonsmooth curve $\hat\Gamma$ as the union $\hat\Gamma = \partial W_\sigma \cup r\cdot \partial W_1$ of the Wulff shape $\partial W_\sigma$ and the circle $ r\cdot \partial W_1$ of a radius $r>0$ touching the Wulff shape from outside at a point $\vecy$ (see Fig~\ref{fig:failurecomparison} (left)). For such a curve we have
\begin{eqnarray*}
\scF_{\hat\Gamma} 
&=& - \frac{L_\sigma(\hat\Gamma)}{2 A(\hat\Gamma)}
= - \frac{L_\sigma(\partial W_\sigma)+L_\sigma(r\cdot\partial W_1)}{2 (A(\partial W_\sigma)+A(r\cdot \partial W_1) )}
\\
&=& - \frac{L_\sigma(\partial W_\sigma)+r L_\sigma(\partial W_1)}{L_\sigma(\partial W_\sigma)+ 2\pi r^2}
\end{eqnarray*}
because $A(\partial W_\sigma) = L_\sigma(\hat\Gamma)/2$. Hence 
$\scF_{\hat\Gamma} > -1 $ provided that the radius $r$ is sufficiently large, $r> L_\sigma(\partial W_1)/2\pi$. For instance, if $\sigma(\nu) = 1 +\varepsilon\cos(m\nu)$, then $r>1$ because $L_\sigma(\partial W_1) =2\pi$ (see section 3). 

Let $\tilde\vecx^0\in\hat\Gamma\cap \partial W_\sigma$ be a point different from $\vecy$ and belonging to a part of the curve representing the Wulff shape. Now, let us construct an initial smooth curve $\tilde\Gamma^0$ which is a continuous perturbation of $\hat\Gamma$, it contains the Wulff shape in the closure of its interior, and such that $\tilde\Gamma^0 \equiv \hat \Gamma$ in some neighborhood ${\mathcal O}(\tilde\vecx)$ of $\tilde\vecx$ (see Fig~\ref{fig:failurecomparison} (right)). The anisoperimetric ratio gradient flow $\tilde\Gamma^t,t\ge 0,$ starting from $\tilde\Gamma^0$ intersects the stationary Wulff shape $\partial W_\sigma$ in the neighborhood ${\mathcal O}(\tilde\vecx)$ for any time $0<t\ll 1$ ($\tilde\Gamma^t$ is plotted by a dashed curve in Fig~\ref{fig:failurecomparison} (right)). This is a consequence of the fact that the normal velocity $\beta$ at $\tilde\vecx^0$ is strictly positive for $t=0$ because $\beta = k_\sigma +\scF_{\tilde\Gamma^0} = 1 +\scF_{\tilde\Gamma^0} >0$ at $\tilde\vecx^0$, and $\tilde\Gamma^0 \cap  {\mathcal O}(\tilde\vecx) = \partial W_\sigma \cap  {\mathcal O}(\tilde\vecx)$, i.e. $k_\sigma =1$ at $\vecx^0$.

\section{Numerical experiments}

In this section we present several examples of evolution of plane curves  minimizing their anisoperimetric ratio. Our scheme belongs to a class of boundary tracking methods taking into account tangential redistribution. In construction of the scheme we employed flowing finite volume discretization in space with a nontrivial tangential velocity, and semi-implicit discretization in time. The advantage of applying a nontrivial tangential velocity consists in its capability to overcome various numerical instabilities of a flow of plane curves like swallow tails and/or merging of numerical grid points leading to break-up of the numerical scheme known for the case when the numerical scheme is constructed with no tangential redistribution. One can find recent progress in \cite{SY2008, SY2011}. Such a scheme is simple and fast, but even if the original problem has a variational structure, it is unclear that the discretized problem has variational structure. On the other hand, in \cite{BenesKimuraYazaki2009}, the authors proposed a semi-discrete scheme with variational structure in a discrete sense. Their scheme has the second order of accuracy in time. However, discretized polygonal curves are restricted to a certain class of curves which is analogous to the admissible class 
in crystalline curvature flows or crystalline algorithm.  In what follows, we propose a hybrid scheme taking into account advantages from both aforementioned  schemes. 

{\bf Discretization scheme.}\ 
For a given initial $N$-sided polygonal curve $\calP^0=\bigcup_{i=1}^N\calS_i^0$, 
we will find a family of $N$-sided polygonal curves $\{\calP^j\}_{j=1, 2, \cdots}$, $\calP^j=\bigcup_{i=1}^N\calS_i^j$, where $\calS_i^j=[\vecx_{i-1}^j, \vecx_{i}^j]$ is the $i$-th edge with $\vecx_{0}^j=\vecx_N^j$ for $j=0, 1, 2, \cdots$. The initial polygon $\calP^0$ is an approximation of $\Gamma^0$ satisfying $\{\vecx_i^0\}_{i=1}^N\subset\calP^0\cap\Gamma^0$, 
and $\calP^j$ is an approximation of $\Gamma^t$ at the time $t=t_j$, 
where $t_j=j\tau$ is the $j$-th discrete time ($j=0, 1, 2, \cdots$) if we use a fixed time increment $\tau>0$, or 
$t_j=\sum_{l=0}^{j-1}\tau_l$ ($j=1, 2, \cdots$; $t_0=0$) if we use 
adaptive time increments $\tau_l>0, l=0, \cdots, j-1 $. 
The updated curve $\calP^{j+1}$ is determined from the data for $\calP^j$ at the previous time step by using discretization in space and time. 
Our two steps scheme will be constructed as follows: in the first step we construct moving polygonal curves which is continuous in time and discrete in space. In the second step we make use of the semi-implicit time discretization scheme for moving polygonal curves. 

{\bf Step 1: Moving polygonal curves.}\ 
Let $\calP(t)=\bigcup_{i=1}^N\calS_i(t)$ be an $N$-sided polygonal curve continuously in time with $\calP(0)=\calP^0$, where 
$\calS_i(t)=[\vecx_{i-1}(t), \vecx_{i}(t)]$ is the $i$-th edge and $\vecx_i(t)$ is the $i$-th vertex 
($i=1, 2, \cdots, N$; $\vecx_0(t)=\vecx_N(t)$). 
The length of $\calS_i$ is denoted by $r_i=|\vecx_{i}-\vecx_{i-1}|$. 
The $i$-th unit tangent vector $\vecT_i$ can be defined as $\vecT_i=(\vecx_{i}-\vecx_{i-1})/r_i$, and 
the $i$-th unit inward normal vector $\vecN_i=\vecT_i^{\bot}$, where $(a, b)^{\bot}=(-b, a)$. 
Then the $i$-th unit tangent angle $\nu_i$ is obtained from 
$\vecT_i=(\cos\nu_i, \sin\nu_i)^{\mathrm{T}}$ in the following way: 
Firstly, from $\vecT_1=(T_{11}, T_{12})^{\mathrm{T}}$, 
we obtain $\nu_1=-\arccos(T_{11})$ if $T_{12}<0$; $\nu_1=\arccos(T_{11})$ if $T_{12}\geq 0$. 
Secondly, for $i=1, 2, \cdots, N$  we successively compute $\nu_{i+1}$ from $\nu_{i}$: 
\[
\nu_{i+1}=\left\{\begin{array}{@{}ll}
\nu_i-\arccos(I), & \mbox{if $D<0$}, \\
\nu_i+\arccos(I), & \mbox{if $D>0$}, \\
\nu_i, & \mbox{otherwise}, 
\end{array}\right.
\]
where $D=\det(\vecT_i, \vecT_{i+1}), 
I=\vecT_i\cdot\vecT_{i+1}$. Finally, we obtain 
$\nu_0=\nu_1-(\nu_{N+1}-\nu_{N})$. 
Then the $i$-th unit inward normal vector $\vecN_i$ is $\vecN_i=(-\sin\nu_i, \cos\nu_i)^{\mathrm{T}}$. 

Let us introduce the ``dual'' volume 
$\calS_i^*=[\vecx_{i}^*, \vecx_{i}]\cup[\vecx_{i}, \vecx_{i+1}^*]$ of $\calS_i$, 
where $\vecx_{i}^*=(\vecx_{i-1}+\vecx_{i})/2$ is the mid point of the $i$-th edge $\calS_i$ 
($i=1, 2, \cdots, N$; $\vecx_{N+1}^*=\vecx_1^*$). 
The length of $\calS_i^*$ is $r_i^*=(r_i+r_{i+1})/2$. 
Then the total length of $\calP$ is 
$L=\sum_{i=1}^Nr_i=\sum_{i=1}^Nr_i^*$, and 
the enclosed area of $\calP$ is 
$A=-\sum_{i=1}^N(\vecx_i\cdot\vecN_i)r_i/2=\sum_{i=1}^N\vecx_{i-1}^{\bot}\cdot\vecx_i/2$. 

We define the $i$-th unit tangent angle of $\calS_i^*$ by $\nu_i^*=(\nu_i+\nu_{i+1})/2=\nu_i+\phi_i/2$, where 
$\phi_i=\nu_{i+1}-\nu_i$ is the angle between the adjacent two edges. 
Then the $i$-th tangent vector at the vertex $\vecx_i$ is 
$\vecT_i^*=(\cos\nu_i^*, \sin\nu_i^*)^{\mathrm{T}}$ and the inward normal vector 
$\vecN_i^*=(-\sin\nu_i^*, \cos\nu_i^*)^{\mathrm{T}}$. 
Hereafter we will use the following abbreviations: 
\[
c_i = \cos\frac{\phi_i}{2}, \quad
s_i = \sin\frac{\phi_i}{2} \quad (i=1, 2, \cdots, N). 
\]
Then it is easy to check that 
\[
\vecT_i^*=c_i\vecT_i+s_i\vecN_i, \quad
\vecN_i^*=c_i\vecN_i-s_i\vecT_i \quad (i=1, 2, \cdots, N). 
\]

The evolution equations of $\calP(t)$ read as follows:
\begin{equation} \label{eq:dot{vecx}_i}
\dot{\vecx}_i=\alpha_i\vecT_i^*+\beta_i\vecN_i^* \quad (i=1, 2, \cdots, N), 
\end{equation}
where $\alpha_i$ and $\beta_i$ are quantities defined on $\calS_i^*$. 
Here and hereafter, we denote $\dot{u}=\du/\dt$. 
The tangential velocities $\{\alpha_i\}$ are defined below and 
the $i$-th normal velocity $\beta_i$ is defined such as
\begin{equation} \label{eq:beta_i}
\beta_i=\frac{\beta_i^*+\beta_{i+1}^*}{2c_i} \quad (i=1, 2, \cdots, N), 
\end{equation}
where the $i$-th normal velocity $\beta_i^*$ is defined on $\calS_i$. It is an approximation of (\ref{geomrovnonloc}) such as
\[
\beta_i^*=\delta_ik_i + \scF_\Gamma \quad (i=1, 2, \cdots, N). 
\]
Here $k_i$ is the $i$-th curvature and the constant value on $\calS_i$ defined as
\[
k_i=\frac{\tan(\phi_i/2)+\tan(\phi_{i-1}/2)}{r_i} \quad (i=1, 2, \cdots, N), 
\]
which is the same as the polygonal curvature in \cite{BenesKimuraYazaki2009}, and 
$\delta_i$ is an approximation of $\delta(\nu_i)$ defined later. 
Then we obtain the time evolution of the total length of $\calP(t)$:
\[
\dot{L}
=-2\sum_{i=1}^N\beta_is_i
=-\sum_{i=1}^Nk_i\beta_i^*r_i, 
\]
and the time evolution of the enclosed area of $\calP(t)$:
\begin{align}
\dot{A}
&=-\sum_{i=1}^N\beta_ic_ir_i^*+\sum_{i=1}^N\alpha_is_i\frac{r_{i+1}-r_i}{2} 
\\
&=-\sum_{i=1}^N\beta_i^*r_i+\mathrm{err}_A ,
\label{eq:semi-discrete-dA/dt} \\
\mathrm{err}_A
&= -\sum_{i=1}^N\beta_i^*\frac{r_{i+1}-2r_i+r_{i-1}}{4}+\sum_{i=1}^N\alpha_is_i\frac{r_{i+1}-r_i}{2}.
\nonumber
\end{align}
These identities represent a discrete version of equations (\ref{eq:area}) provided that the distribution $r_i\equiv L/N, i=1, 2, \cdots, N,$ is uniform because the error term $\mathrm{err}_A=0$ is vanishing.

To realize this uniform distribution asymptotically, we assume that
\[
r_i-\frac{L}{N}=\eta_ie^{-f(t)} 
\]
\[
\left(\mbox{$\sum_{i=1}^N\eta_i=0$, $f(t)\to\infty$ as $t\to T_{max}\leq \infty$}\right). 
\]
By using a relaxation term $\omega(t)=f'(t)$ we obtain
\begin{equation}
\dot{r}_i-\frac{\dot{L}}{N}
=\left(\frac{L}{N}-r_i\right)\omega(t), \quad \int_0^{T_{max}}\omega(t)\,\dt=\infty 
\label{eq:relaxation-omega}
\end{equation}
$(i=1, 2, \cdots, N)$. Taking into account the relations:
\begin{align*}
\dot{r}_i
&= (\dot{\vecx}_i-\dot{\vecx}_{i-1})\cdot\vecT_i
\\
&= -\beta_is_i-\beta_{i-1}s_{i-1}+c_i\alpha_i-c_{i-1}\alpha_{i-1}
\\
&= \frac{\dot{L}}{N}+\left(\frac{L}{N}-r_i\right)\omega(t), 
\end{align*}
we deduce $N-1$ equations for tangential velocities $\alpha_i$ ($i=2, 3, \cdots, N$):
\begin{align*}
& \alpha_i=\frac{\Psi_i}{c_i}+\frac{c_1}{c_i}\alpha_1 \quad (i=2, 3, \cdots, N), \\
& \Psi_i=\psi_2+\psi_3+\cdots+\psi_i \quad (i=2, 3, \cdots, N), \\
& \psi_i=\beta_is_i+\beta_{i-1}s_{i-1}-\frac{2}{N}\sum_{i=1}^N\beta_is_i+\left(\frac{L}{N}-r_i\right)\omega(t). 
\end{align*}
To determine $\alpha_1$, we add one more linear equation of the form $\sum_{i=1}^N\alpha_ip_i=P$, 
which is independent of the above $N-1$ equations. Since
\begin{equation} \label{eq:RQ}
R=c_1\sum_{i=1}^N\frac{p_i}{c_i}, \quad Q=\sum_{i=2}^N\frac{p_i}{c_i}\Psi_i, 
\end{equation}
we obtain $\alpha_1=(P-Q)/R$. 
Next, we propose three candidates for each $\{p_i\}$ and $P$, and choose one of them in the following way: 

{\it Candidate 1.}\ We put 
\[
p_i=s_i\frac{r_{i+1}-r_i}{2} \quad (i=1, 2, \cdots, N), 
\]
\[
P=\sum_{i=1}^N\beta_i^*\frac{r_{i+1}-2r_i+r_{i-1}}{4}, 
\]
and from (\ref{eq:RQ}) we calculate $R$ and $Q$. We denote this $R$ by $R_1$. 
If the above equation holds, then $\mathrm{err}_A=0$ and $\dot{A}=-\sum_{i=1}^N\beta_i^*r_i$ hold in (\ref{eq:semi-discrete-dA/dt}). 
However, if distribution of grid points are almost uniform, then the above equation is almost nothing. 
Therefore we need another candidate. 

{\it Candidate 2.}\ 
For the $i$-th quantities ${\sf F}_i$ defined on $\calS_i$ and ${\sf G}_i$ defined on $\calS_i^*$, 
we define the average along $\calP$ such as
\[
\langle{\sf F}\rangle = \frac{1}{L}\sum_{i=1}^N{\sf F}_ir_i, \quad
\langle{\sf G}\rangle^* = \frac{1}{L}\sum_{i=1}^N{\sf G}_ir_i^*. 
\]
Since $L=\sum_{i=1}^Nr_i=\sum_{i=1}^Nr_i^*$, we have $\langle 1\rangle=\langle 1\rangle^*=1$. 
Moreover, for $\alpha_i^*=(\alpha_i+\alpha_{i-1})/2$ defined on $\calS_i$, 
the relation $\langle \alpha\rangle^*=\langle \alpha^*\rangle$ holds. 

The second candidate of linear equation is the zero-average $\langle \alpha \rangle^*=0$, 
that is, $p_i=r_i^*$ for $i=1, 2, \cdots, N$ and $P=0$. 
From this and (\ref{eq:RQ}) we calculate $R$ and $Q$. We denote this $R$ by $R_2$. 

The purpose of this section is to present numerical simulations of the geometric flow evolving according to the evolution equation (\ref{eq:isoperimnormalaniso}). 
Before we introduce the third candidate, 
we calculate the discrete version of (\ref{eq:isoperimnormalaniso}). 
Let the total interfacial energy be 
\[
L_\sigma(t)=\sum_{i=1}^N\sigma(\nu_i(t))r_i(t). 
\]
The time derivative of $\vecT_i=(\cos\nu_i, \sin\nu_i)^{\mathrm{T}}$ is 
$\dot{\vecT}_i=\dot{\nu}_i\vecN_i$. 
Then we have $r_i\dot{\nu}_i=(\dot{\vecx}_i-\dot{\vecx}_{i-1})\cdot\vecN_i$, 
from which it follows that
\begin{align*}
\dot{L}_\sigma
&= \sum_{i=1}^N(\sigma'_i\dot{\nu}_ir_i+\sigma_i\dot{r}_i)
= -\sum_{i=1}^N{k_{\sigma}}_i\beta_i^*r_i+\sum_{i=1}^N{k_{\sigma}}_i\tilde{p}_i\alpha_i, \\
\tilde{p}_i
&=(\sigma'_i+\sigma'_{i+1})s_i+(\sigma_i-\sigma_{i+1})c_i. 
\end{align*}
Here $\sigma_i=\sigma(\nu_i)$, $\sigma'_i=\sigma'(\nu_i)$, and 
${k_{\sigma}}_i$ is discrete version of the weighted curvature in the following sense: 
\[
{k_{\sigma}}_i
=\delta_ik_i,
\]
\[
\delta_i=\frac{\sigma'_{i+1}-\sigma'_{i-1}}{2(t_i+t_{i-1})}
+\frac{\sigma_{i+1}t_i
+\sigma_{i}(t_i+t_{i-1})
+\sigma_{i-1}t_{i-1}}{2(t_i+t_{i-1})},
\]
\[
t_i=\tan\frac{\phi_i}{2}, 
\]
and $\delta_i$ is discrete weight of $\delta(\nu)=\sigma''(\nu)+\sigma(\nu)$ at $\nu=\nu_i$. 
Note that $\delta_i\to\delta(\nu_i)$ holds as $\phi_i, \phi_{i-1}\to0$ formally, and 
even the case where $t_i+t_{i-1}=0$, 
${k_{\sigma}}_i$ is well-defined, 
since $k_i=(t_i+t_{i-1})/r_i$ and then denominator of ${k_{\sigma}}_i$ is $2r_i$. 

We obtain
\begin{align*}
\frac{\mathrm{d}}{\dt}\frac{L_\sigma^2}{A}
&= \frac{2L_\sigma}{A}\left(\dot{L}_\sigma-\frac{L_\sigma}{2A}\dot{A}\right), 
\\
\dot{L}_\sigma-\frac{L_\sigma}{2A}\dot{A}
&= -\sum_{i=1}^N\left({k_{\sigma}}_i-\frac{L_\sigma}{2A}\right)\beta_i^*r_i
+\mathrm{err}_{ratio}, \\
\mathrm{err}_{ratio}
&= \sum_{i=1}^N\tilde{p}_i\alpha_i
\\
&+\frac{L_\sigma}{2A}\biggl(
\sum_{i=1}^N\beta_i^*\frac{r_{i+1}-2r_i+r_{i-1}}{4}
\\
&-\sum_{i=1}^N\alpha_is_i\frac{r_{i+1}-r_i}{2}
\biggr). 
\end{align*}

{\it Candidate 3.}\ We put 
\[
p_i=s_i\frac{r_{i+1}-r_i}{2}-\frac{2A}{L_\sigma}\tilde{p}_i \quad (i=1, 2, \cdots, N), 
\]
\[
P=\sum_{i=1}^N\beta_i^*\frac{r_{i+1}-2r_i+r_{i-1}}{4}, 
\]
and from (\ref{eq:RQ}) we calculate $R$ and $Q$. We denote this $R$ by $R_3$. 
Let the $i$-th normal velocity defined on $\calS_i$ be
\[
\beta_i^*={k_{\sigma}}_i-\frac{L_\sigma}{2A} \quad (i=1, 2, \cdots, N), 
\]
which is discrete version of (\ref{eq:isoperimnormalaniso}). 
If we choose the candidate 3 and its hold exactly, 
then $\mathrm{err}_{ratio}=0$ holds and we obtain
\[
\frac{\mathrm{d}}{\dt}\frac{L_\sigma^2}{A}
=-\frac{2L_\sigma}{A}
\sum_{i=1}^N\left({k_{\sigma}}_i-\frac{L_\sigma}{2A}\right)^2r_i<0. 
\]

{\it Choice one from three candidates.}\ 
We choose candidate number $l$ satisfying 
\[
|R_l|=\max\{|R_1|, |R_2|, |R_3|\}. 
\]

\begin{remark}
If, for definition of the normal velocity defined on $\calS_i^*$, we use
\[
\beta_i=\frac{\beta_i^*r_i+\beta_{i+1}^*r_{i+1}}{2c_ir_i^*}
\]
instead of (\ref{eq:beta_i}), 
then ${k_{\sigma}}_i$ can not be divided into weighted part $w_i$ and the curvature $k_i$. 
Moreover, if we use 
\[
\beta_i=\frac{\beta_i^*+\beta_{i+1}^*}{2}
\]
instead of (\ref{eq:beta_i}), 
then ${k_{\sigma}}_i$ can be divided into weighted part and the curvature. 
However, we have $\mathrm{err}_A \neq 0$, that is $\dot{A}=-\sum_{i=1}^N\beta_i^*r_i$ does not hold 
even if uniform distribution $r_i\equiv L/N$ holds for all $i=1, 2, \cdots, N$. 
\end{remark}

{\bf Step 2: Discretization in time.}\ 
We use semi-implicit scheme for discretization of (\ref{eq:dot{vecx}_i}). 
Next we develop  expression (\ref{eq:dot{vecx}_i}) as follows: 
\begin{eqnarray}
\dot{\vecx}_i
&=& \alpha_i\vecT_i^*+\beta_i\vecN_i^* \nonumber \\
&=& \frac{1}{2}\left(\frac{\alpha_i}{c_i}-\beta_is_i\right)\vecT_{i}
+ \frac{1}{2}\left(\frac{\alpha_i}{c_i}+\beta_is_i\right)\vecT_{i+1} 
\nonumber \\
&&+ \beta_ic_i\frac{\vecN_{i}+\vecN_{i+1}}{2} \label{eq:dot{vecx}3} 
 \label{eq:dot{vecx}4}\\
&=& \frac{1}{2}\left(\frac{\alpha_i}{c_i}-\frac{\beta_i}{s_i}\right)\vecT_{i}
+ \frac{1}{2}\left(\frac{\alpha_i}{c_i}+\frac{\beta_i}{s_i}\right)\vecT_{i+1}. \nonumber
\end{eqnarray}
Here we have used the relation
\begin{eqnarray*}
\vecT_i^*
&=&\frac{\vecT_{i+1}+\vecT_{i}}{2c_i}, \quad
\\
\vecN_i^* 
&=&s_i\frac{\vecT_{i+1}-\vecT_{i}}{2}
+c_i\frac{\vecN_{i+1}+\vecN_{i}}{2}
=\frac{\vecT_{i+1}-\vecT_{i}}{2s_i}. 
\end{eqnarray*}
Let $\mu$ be a parameter satisfying $\mu=0$ if $\min_{1\leq i\leq N}|s_i|=0$, otherwise $\mu\in (0, 1]$. 
For the parameter $\mu\in[0, 1]$ we use the linear interpolation of (\ref{eq:dot{vecx}3}) and (\ref{eq:dot{vecx}4}), 
since we can not use (\ref{eq:dot{vecx}4}) if $s_i=0$. 

Put $b_i=(1-\mu)s_i+\mu/s_i$ for $\mu\in (0, 1]$ and $b_i=s_i$ for $\mu=0$. 
Then the evolution equation instead of (\ref{eq:dot{vecx}_i}) will be
\begin{eqnarray}
\dot{\vecx}_i
&=&\frac{1}{2}\left(\frac{\alpha_i}{c_i}-\beta_ib_i\right)\vecT_{i}
+\frac{1}{2}\left(\frac{\alpha_i}{c_i}+\beta_ib_i\right)\vecT_{i+1}
\nonumber \\
&& +\beta_ic_i(1-\mu)\frac{\vecN_{i}+\vecN_{i+1}}{2}. 
\label{eq:dot{vecx}5}
\end{eqnarray}
From $\vecT_i=(\vecx_i-\vecx_{i-1})/r_i$, 
we discretize (\ref{eq:dot{vecx}5}) in time and obtain the following tridiagonal linear system under 
the periodic boundary condition:
\begin{eqnarray*}
\frac{\vecx_i^{j+1}-\vecx_i^{j}}{\tau_j}
&=&-a_-\vecx_{i-1}^{j+1}+a_0\vecx_i^{j+1}-a_+\vecx_{i+1}^{j+1}
\\
&&+\beta_i^jc_i^j(1-\mu)\frac{\vecN_{i}^j+\vecN_{i+1}^j}{2}. 
\end{eqnarray*}
Here $a_0=a_-+a_+$ and 
\[
a_-=\frac{1}{2r_i^j}\left(\frac{\alpha_i^j}{c_i^j}-\beta_i^jb_i^j\right),\]
\[
a_+=-\frac{1}{2r_{i+1}^j}\left(\frac{\alpha_i^j}{c_i^j}+\beta_i^jb_i^j\right) 
\]
for $i=1, 2, \cdots, N$ and $j=0, 1, 2, \cdots$. 
For the choice of $\mu$ we use
\[
\mu=\frac{\min_i|s_i^j|}{\max_i|s_i^j|}\in [0, 1]. 
\]
Here and hereafter, 
$\min_i$ and $\max_i$ mean $\min_{1\leq i\leq N}$ and $\max_{1\leq i\leq N}$, respectively. 
Note that $\max_i|s_i^j|>0$ holds for closed curves. 

In order to ensure solvability of the above linear system, we require a simple condition on the diagonal dominance. Adopting such a condition the adaptive time step $\tau_j$ satisfies
\[
\tau_j=\frac{\min_i r_i^j}
{2(1+\lambda)(\max_i|\alpha_i^j/c_i^j|+\max_i|\beta_i^jb_i^j|)} \quad
(\lambda>0). 
\]

\begin{figure}[ht]
\begin{center}
\begin{tabular}{@{}c@{}}
\scalebox{0.65}{\includegraphics{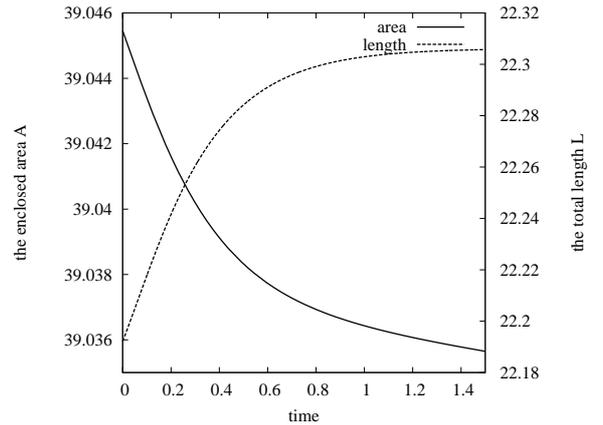}}
\end{tabular}
\end{center}
\caption{%
Initial decrease of the enclosed area $A$ and increase of the total length $L$. 
}
\label{fig:Wulff-sigma-bar-peak6-AL-vs-time}
\end{figure}

{\bf Simulation.}\ 
In following all figures, for the prescribed $\widehat{\tau}>0$, we plot every $\mu\widehat{\tau}$ discrete 
time step using discrete points representing the evolving curve. 
In every $3\mu\widehat{\tau}$ time step, we plot a polygonal curve connecting those points, 
where $\mu=[[T/\widehat{\tau}]/100]$ ($[x]$ is the integer part of $x$), 
and $T=1.5$ is the final computational time. 
We use $\omega\equiv 1000$ as the relaxation term in (\ref{eq:relaxation-omega}). 
Note that if we use small $\omega$, 
then asymptotic speed for uniform distribution becomes slow, and 
for some choice of $\sigma$ area-decreasing phenomena do not hold near the initial time (cf. Fig~\ref{fig:Wulff-sigma-bar-peak6-AL-vs-time}~(a)). 

\begin{figure}[ht]
\begin{center}
\includegraphics[width=1.5in]{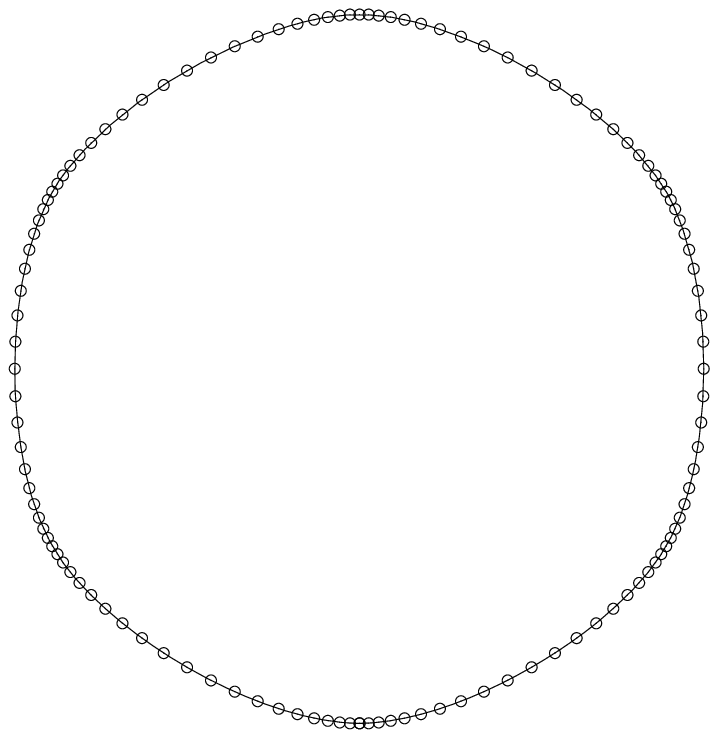}
\quad
\includegraphics[width=1.5in]{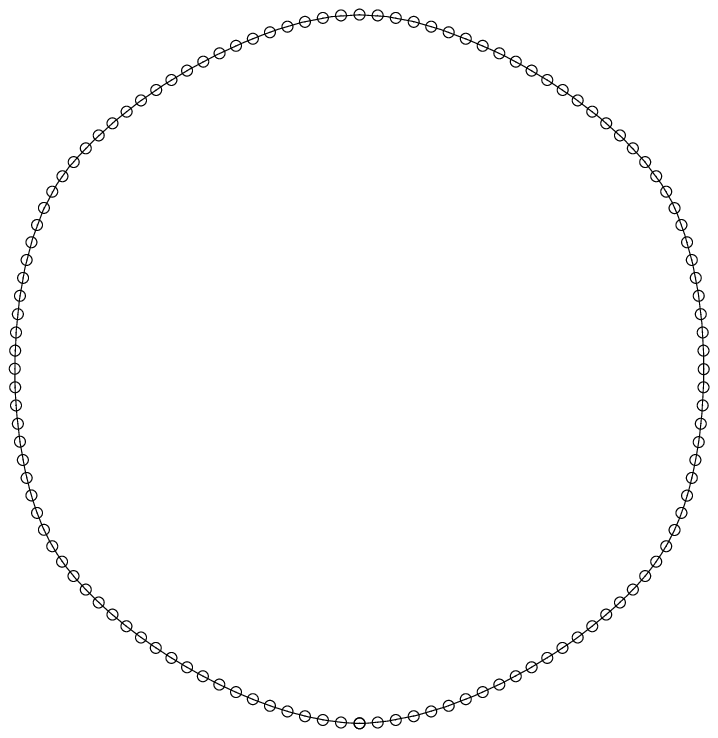}
\\
(a) \hskip 1in (b) 
\\
\includegraphics[width=1.5in]{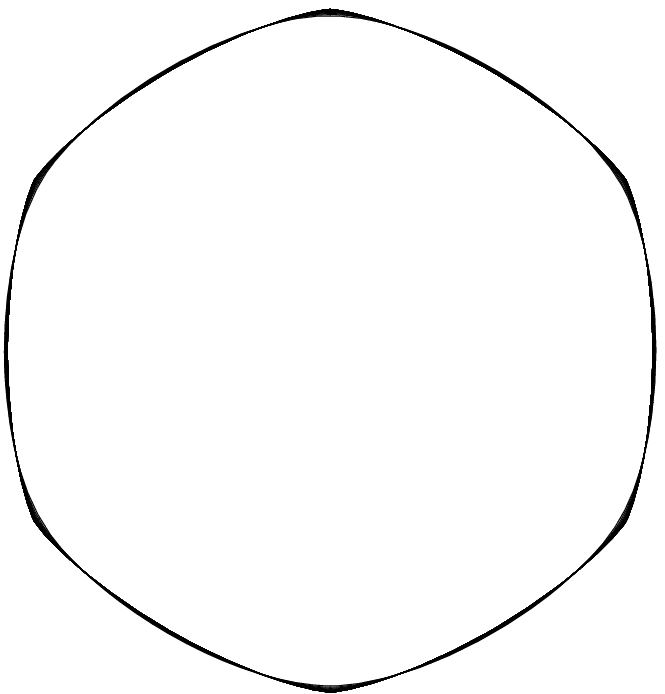}
\\
(c)
\end{center}
\caption{%
(a) The Wulff shape $\partial W_{\bar{\sigma}}$ with $N=120$ points, 
(b) its uniform parameterization and its time evolution (c). 
}
\label{fig:Wulff-sigma-bar-peak6}
\end{figure}

{\bf Wulff shapes and area-decreasing phenomenon.}\ 
If  $\sigma(\nu) =1 + \varepsilon \cos(m \nu)$ is the anisotropy density function of the degree $m$ then, by using (\ref{eq:areaW}), we obtain the explicit expression for the mixed anisotropy function:
\[
\bar{\sigma} = \sqrt{\pi}\, \sigma + \sqrt{|W_\sigma|}
\]
\[
=\sqrt{\pi}\left(1+\eps\cos m\nu+\sqrt{1-\eps^2\frac{m^2-1}{2}}\right). 
\]
In order to verify the area-decreasing and the length-increasing phenomenon at the initial time as in Theorem~\ref{th:areadecrasing}, 
we use $\partial W_{\bar{\sigma}}$ as the initial curve. 
Fig~\ref{fig:Wulff-sigma-bar-peak6}~(a) indicates $\partial W_{\bar{\sigma}}$ with $m=6$. 
Its discretization is given by the uniform $N$-division of the $u$-range $[0, 1]$. 
Fig~\ref{fig:Wulff-sigma-bar-peak6}~(b) indicates the same Wulff shape, but the grid points are distributed uniformly. 
Fig~\ref{fig:Wulff-sigma-bar-peak6}~(c) indicates the time evolution starting from (b). 

\begin{figure}[ht]
\begin{center}
\includegraphics[width=1.65in]{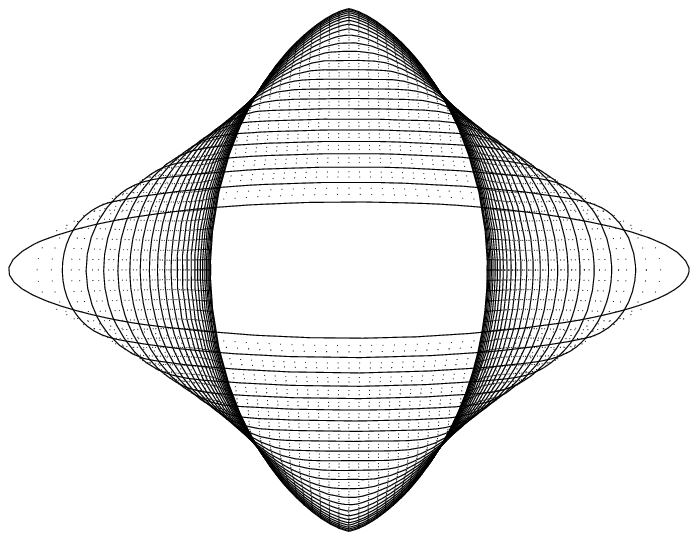}
\\
(a)
\\ \ 
\\
\includegraphics[width=1.65in]{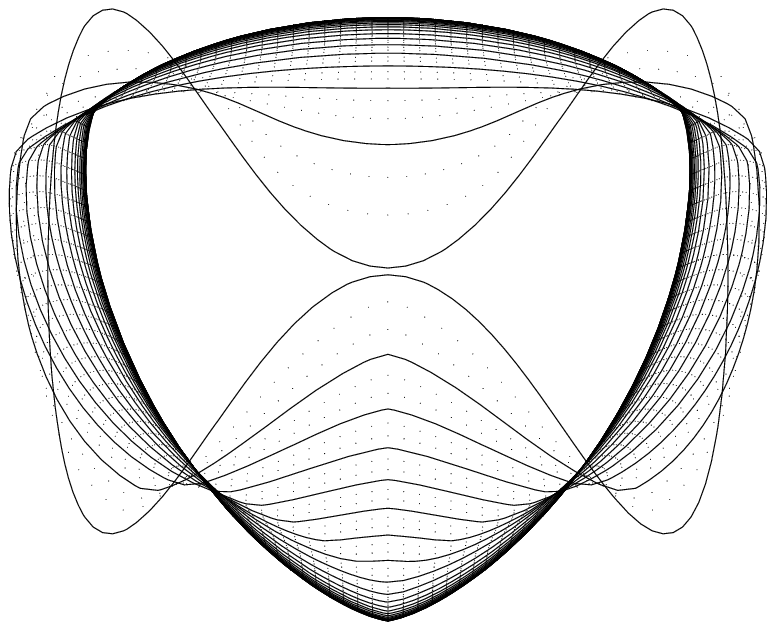}
\\
(b) 
\\ \ 
\\
\includegraphics[width=1.65in]{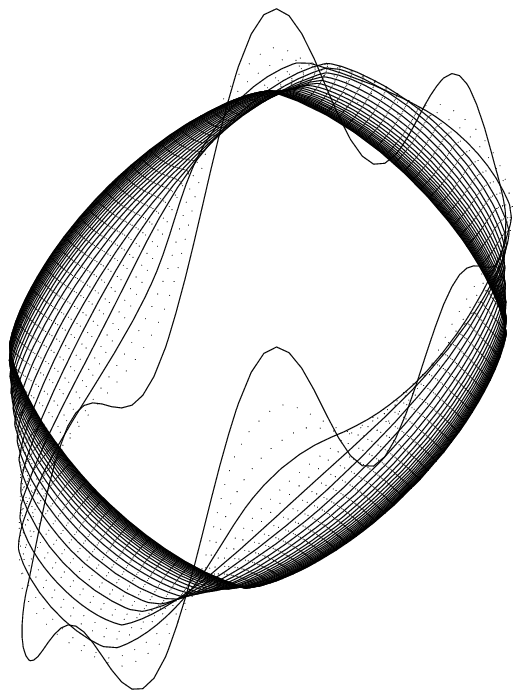}
\\
(c)
\\ \ 
\\
\includegraphics[width=1.65in]{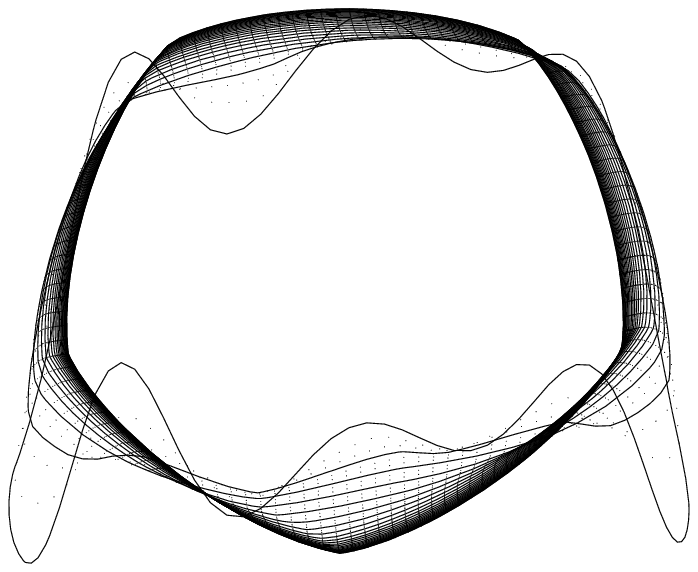} 
\\
(d)
\\ \ 
\\
\includegraphics[width=1.65in]{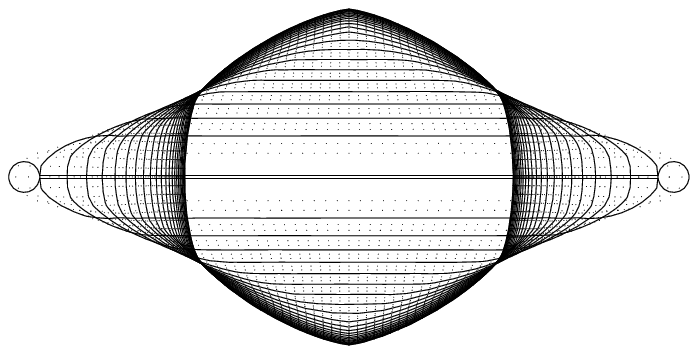}
\\
(e)
\\
\end{center}
\caption{%
Evolution of curves starting from the initial curves with various choice of peak of $\sigma$. 
}
\label{fig:evolution-of-curves}
\end{figure}

Although deformation of $\partial W_{\bar{\sigma}}$ is very small (see Fig~\ref{fig:Wulff-sigma-bar-peak6}~(c)), the area-decreasing and the length-increasing phenomenon can numerically verified by using the aforementioned numerical discretization scheme. The behavior of the enclosed area and total length of curves evolved from the initial Wulff shape $\partial W_{\bar\sigma}$ with the mixed anisotropy density function is shown in Fig~\ref{fig:Wulff-sigma-bar-peak6-AL-vs-time}. 

{\bf Initial test curves.}\ 
As initial test examples we use the boundary $\partial W_{\bar{\sigma}}$ of the Wulff shape as well as
the following initial curves $\vecx(u, 0)=(x_1(u), x_2(u))^{\mathrm{T}}$ ($u\in[0, 1]$) parameterized by
\begin{align*}
&\mathrm{(a)}&
& \mbox{ellipse}: 
x_1(u)=a\cos 2\pi u, \ 
x_2(u)=b\sin 2\pi u, \\
&\mathrm{(b)}&
& \mbox{dumbbell shape}: z = 2\pi u,  \\
& & &
x_1(u)=\cos z, \ 
x_2(u)=2.0\sin z - 1.99\sin^3z, \ 
\\
&\mathrm{(c)}&\ 
& x_1(u)=\cos z, \ x_2(u)=0.7\sin z+\sin x_1+x_3^2, \ 
\\
& & &
x_3=\sin(3z)\sin z, \ z=2\pi u, \\
&\mathrm{(d)}&
& 
x_1(u)=1.5\cos z, \ 
\\
& & &
x_2(u)=1.5(0.6\sin z+0.5x_3^2+0.4\sin x_4+0.1\sin x_5), \\
&&&
x_3=\sin(3z)\sin z, \ x_4=2x_1^2, \ x_5=3e^{-x_1}, \ z=2\pi u, \\
&\mathrm{(e)}&
& \mbox{thin-dumbbell shape}: \\
&&&
(x_1(u), x_2(u))=\left\{\begin{array}{@{}ll}
(\bar{x}_1(u), \bar{x}_2(u)) & \\ \qquad \hbox{for}\  0\leq u<0.5 & \\
-(\bar{x}_1(u-0.5), \bar{x}_2(u-0.5)) \\ \qquad \hbox{for}\  0.5\leq u\leq 1 &
\end{array}\right., \\
&&&
(\bar{x}_1(u), \bar{x}_2(u))=\left\{\begin{array}{@{}ll}
(\hat{x}_1(u), \hat{x}_2(u)) \\ \qquad \hbox{for}\ 0\leq u<0.25 & \\
(-\hat{x}_1(0.5-u), \hat{x}_2(0.5-u)) \\ \qquad \hbox{for}\ 0.25\leq u\leq 0.5 &
\end{array}\right., \\ 
&&&
\hat{x}_1(u)=\left\{\begin{array}{@{}ll}
beam+rad(1+\cos(8(\pi-\theta_\eps)u)) \\ \quad  \hbox{for}\ 0\leq u<0.125& \\
2(beam+rad(1-\cos\theta_\eps))(1-4u) \\ \qquad  \hbox{for}\ 0.125\leq u\leq 0.25 &
\end{array}\right., \\
&&&
\hat{x}_2(u)=\left\{\begin{array}{@{}ll}
rad\sin(8(\pi-\theta_\eps)u)) \\ \qquad  \hbox{for}\ 0\leq u<0.125 & \\
\eps \\ \qquad  \hbox{for}\ 0.125\leq u\leq 0.25 &
\end{array}\right., 
\end{align*}
where $beam>0$ and $0<\eps<rad$ are parameters and $\theta_\eps=\arcsin(\eps/rad)$. 

In all examples, the initial discretization is given by the uniform $N$-division of the $u$-range $[0, 1]$. Fig~\ref{fig:evolution-of-curves} indicates numerical simulation with the initial curves. We choose several peaks of $\sigma$ such as  (a) $m=2$, (b) $m=3$, (c) $m=4$, (d) $m=5$, and (e) $m=6$.

\begin{figure}[ht]
\begin{center}
\includegraphics[width=1.8in]{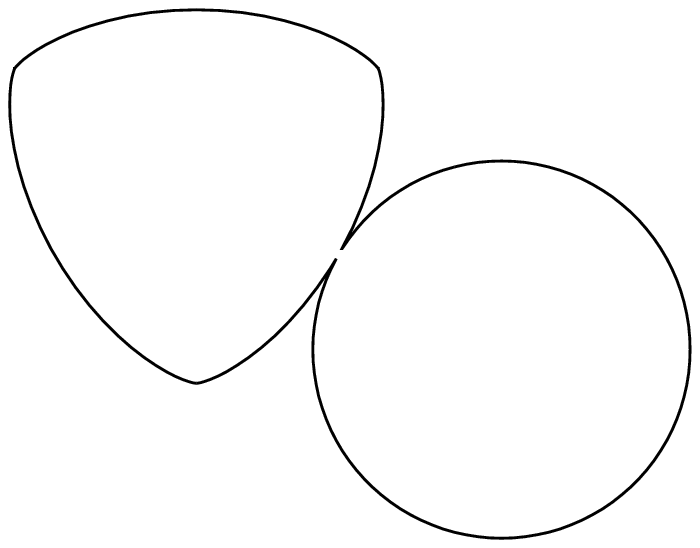}
\\
(a)
\\ \ 
\\
\includegraphics[width=1.8in]{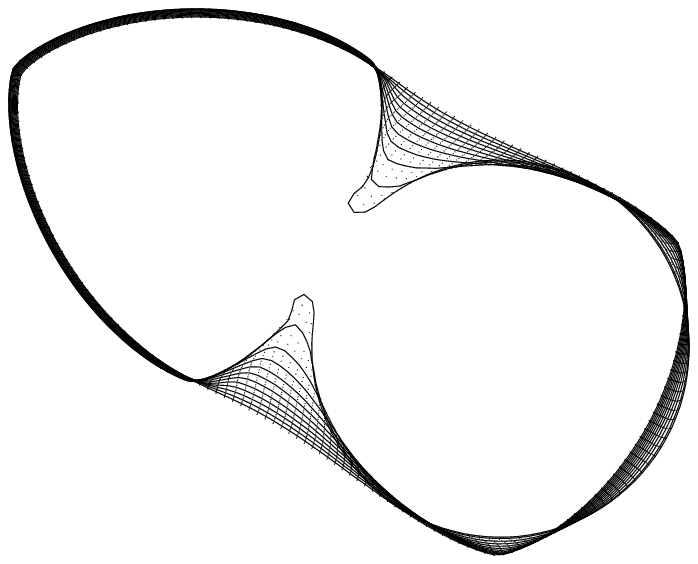}
\\
(b)
\\ \ 
\\
\includegraphics[width=1.8in]{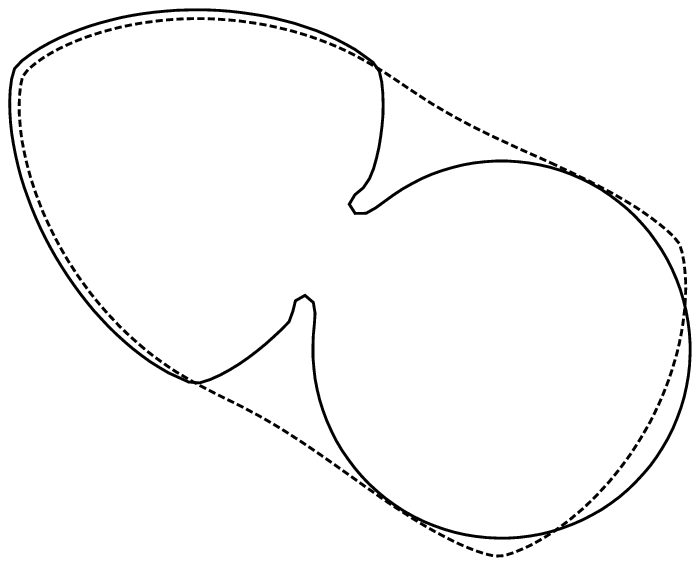}
\\
(c)
\end{center}
\caption{%
Anisoperimetric ratio gradient flow breaking comparison principle. 
}
\label{fig:comparison-breaking}
\end{figure}

{\bf Breaking of a comparison principle.}\ 
In Fig~\ref{fig:comparison-breaking} we plot the initial curve consisting of the union of the boundary of a Wulff shape $\partial W_\sigma$ (with $\sigma$  having degree $m=3$) touched by a circle with a sufficiently large circle with a radius $r>1$ (see section 7). As soon as it evolved by the anisoperimetric ratio gradient flow it intersects the stationary Wulff shape $\partial W_\sigma$. 
Numerically computed examples displayed in  Fig~\ref{fig:comparison-breaking} (a) and (c) correspond to those of the conceptual Fig~\ref{fig:failurecomparison}.

\section{Conclusions}
In this paper, we have derived and analyzed a gradient flow of closed planar curves minimizing the anisoperimetric ratio. A geometric law for normal velocity is a function of the anisotropic curvature and it depends on the total interfacial energy and enclosed area of the curve. We also derived a new mixed anisoperimetric inequality for the product of total interfacial energies corresponding to different anisotropy functions. Interestingly enough, there exist initial curves for which the enclosed area is a decreasing function with respect to time. This is in contrast to the known property of a gradient flow minimizing isoperimetric ratio. We also derived a stable numerical scheme based on the flowing finite volumes method. Theoretical results have been illustrated be several computational examples.


%
\IAENGpeerreviewmaketitle

\ifCLASSOPTIONcaptionsoff
  \newpage
\fi

\end{document}